\documentclass[10pt,a4paper]{article}

 \usepackage{amssymb,amsmath,amsfonts,amssymb, amsthm}
 \usepackage{graphics,graphicx,color}
 -\marginparwidth \oddsidemargin -4mm \evensidemargin
 \oddsidemargin%
 \advance\hoffset by 5mm
 \usepackage{amsmath}    
 \advance\hoffset by 5mm

 \title{Energy Dissipation  and Regularity\\
 for a  Coupled  Navier-Stokes and Q-Tensor  System }
 \date{\today}
\author{Marius Paicu\footnote{Universit\'e Paris-Sud, Laboratoire de Math\'ematiques, 91405 Orsay Cedex, France. E-mail: marius.paicu@math.u-psud.fr }\,\,\,\,  and Arghir Zarnescu \footnote{Mathematical Institute, 24-29 St. Giles', Oxford, OX1 3LB, United Kingdom. E-mail: zarnescu@maths.ox.ac.uk} }

\newtheorem{remark}{Remark}
 \newtheorem{lemma} {Lemma}
 \newtheorem{proposition}{Proposition}
 \newtheorem{theorem} {Theorem}

\begin{document}

  \maketitle
 \begin{abstract} We study a complex non-newtonian fluid that models the flow of nematic liquid crystals. The fluid is described  by a system that couples a forced Navier-Stokes system with a parabolic-type system. We prove the existence of global weak solutions in dimensions two and three. We show the existence of a Lyapunov functional for the smooth solutions of the coupled system and use the cancellations that allow its existence to prove higher global regularity, in dimension two. We also show the weak-strong uniqueness in dimension two. 
\end{abstract}

\section{Introduction}

In this paper we study the global existence of solutions for a system describing the evolution of a nematic liquid crystal flow. The system couples a forced Navier-Stokes system, describing the flow, with a parabolic-type system describing the evolution of the nematic crystal director fields ($Q$-tensors). The coupled system has a Lyapunov functional made of two parts: the free energy due to  the director fields and the kinetic energy of the fluid. This functional describes, from a physical point of view, the dissipation of the energy of the complex fluid.

\par In the first part of the paper we use, in a classical manner, the apriori bounds provided by the energy dissipation to prove the existence of  global weak solutions in the natural energy space. In the second part,   we study the case where the fluid evolves in the two dimensional space and  prove the existence of a global regular solution issued from an appropriately regular initial data. In the two dimensional space we also show that for an appropriately regular initial data the weak and the strong solutions  coincide.
The main contribution of this paper is to show how to use  the specific {\it coupling of the system} (the coupling structure that allows the system to dissipate energy) not only at the level of regularity of weak solutions but also to transport arbitrarily large (enough) regularity of the initial data.  Thus we show that for this type of complex fluids the existence of an energy dissipation is intrinsically related to the high regularity of the solutions.

\par There exist several competing theories that attempt to capture the complexity of nematic liquid crystals, and a comparative discussion and further references are available for instance in \cite{lin&liu}, \cite{mz}. In the present paper we use one of the most comprehensive description of nematics, the $Q$-tensor description, proposed by P.G. de Gennes \cite{dg}. There exist various specific models that all use the $Q$-tensor description and a comparative discussion of the main models is available for instance in \cite{sonnet}.
\par In this paper we  use a model proposed by Beris and Edwards \cite{berised}, that one can find  in the physics literature for instance in \cite{equationsELreduction}, \cite{equations}. An important feature of this model is that if one assumes smooth solutions and one formally takes $Q(x)=s_+(n(x)\otimes n(x)-\frac{1}{3}Id)$, with $s_+$ a constant (depending on the parameters of the system, see for instance \cite{mz})  and $n:\mathbb{R}^d\to\mathbb{S}^{d-1}$ smooth, then the equations reduce (see \cite{equationsELreduction}) to the generally accepted equations of Ericksen, Leslie and Parodi \cite{leslie&ericksen}. The system we study is related structurally to other models of complex fluids coupling a transport equation with a forced Navier-Stokes system  \cite{suli}, \cite{cftz}, \cite{constantin&masmoudi}, \cite{constantin&seregin}, \cite{chemin&masmoudi},  \cite{lions&masmoudi}, \cite{linandco}, \cite{masmoudiFENE}, \cite{schonbek}. In our case the Navier-Stokes equations are coupled with a parabolic type system, but we also have two more derivatives (than in the previously mentioned models)  in the forcing term of the Navier-Stokes equations.  The Ericksen-Leslie-Parodi system describing nematic liquid crystals, whose structure is closer to our system (but  that has one less derivative in the forcing term of the Navier-Stokes equations) was studied in \cite{lin}, \cite{linLE}, \cite{linarma}.

\par In the following we  use a partial Einstein summation convention, that is we  assume summation over repeated {\it greek} indices, but not over the repeated  {\it latin} indices. We consider the equations as described in \cite{equationsELreduction},~\cite{equations} but assume that the fluid has constant  density in time.
\par We denote

\begin{equation}
 S(\nabla u,Q)\stackrel{def}{=}(\xi D+\Omega)(Q+\frac{1}{3}Id)+(Q+\frac{1}{3}Id)(\xi D-\Omega)-2\xi (Q+\frac{1}{3}Id)\textrm{tr}(Q\nabla u)
 \label{rel:defS}
\end{equation} where $D\stackrel{def}{=}\frac{1}{2}\left(\nabla u+(\nabla u)^T\right)$ and $\Omega\stackrel{def}{=}\frac{1}{2}
\left(\nabla u-(\nabla u)^T\right)$ are the symmetric part and the antisymmetric part, respectively, of the velocity gradient tensor $\nabla u$. The term $S(\nabla, Q)$ appears in the equation of motion of the order-parameter, $Q$, and describes how the flow gradient rotates and stretches the order-parameter. The constant $\xi$ depends on the molecular details of a given liquid crystal and measures the ratio between the tumbling and the aligning effect that a shear flow would exert over the liquid crystal directors.

\par We also denote:

\begin{equation}
 H\stackrel{def}{=}-aQ+b[Q^2-\frac{\textrm{tr}(Q^2)}{3}Id]-cQ\textrm{tr}(Q^2)+L\Delta Q
\label{rel:defH}
\end{equation} where $L>0$.

\par With the notations above we have the coupled system:

\begin{equation}
\left\{\begin{array}{l}
      (\partial_t+u\cdot \nabla)Q-S(\nabla,Q)=\Gamma H \\
      \partial_t u_\alpha+u_\beta\partial_\beta u_\alpha=\nu\partial_{\beta\beta} u_\alpha+\partial_\alpha p+\partial_\beta \tau_{\alpha\beta}+\partial_\beta\sigma_{\alpha\beta}\\
       \partial_\gamma u_\gamma=0
     \end{array}\right.
     \label{eq:allsystem}
\end{equation} where $\Gamma>0,\nu>0$ and  we have the symmetric part of the additional stress tensor:

\begin{equation}
\tau_{\alpha\beta}=-\xi \left(Q_{\alpha\gamma}+\frac{\delta_{\alpha\gamma}}{3}\right)H_{\gamma\beta}-\xi H_{\alpha\gamma}\left(Q_{\gamma\beta}+\frac{\delta_{\gamma\beta}}{3}\right)+2\xi (Q_{\alpha\beta}+\frac{\delta_{\alpha\beta}}{3})Q_{\gamma\delta}H_{\gamma\delta}-L\left(\partial_\beta Q_{\gamma\delta}\partial_\alpha Q_{\gamma\delta}+\frac{\delta_{\alpha\beta}}{3}Q_{\nu\varepsilon}Q_{\nu\varepsilon}\right)\end{equation} and an antisymmetric part:

\begin{equation}
\sigma_{\alpha\beta}=Q_{\alpha\gamma}H_{\gamma\beta}-H_{\alpha\gamma}Q_{\gamma\beta}
\end{equation}

\par In the rest of the paper we restrict ourselves to the case $\xi=0$. This means that the molecules are such that they only tumble in a shear flow, but are not aligned by such a flow. In this case the system (\ref{eq:allsystem}) reduces to:

\begin{equation}
\left\{\begin{array}{l}
      (\partial_t+u_\gamma\cdot \partial_\gamma)Q_{\alpha\beta}-\Omega_{\alpha\gamma} Q_{\gamma\beta}+Q_{\alpha\gamma}\Omega_{\gamma\beta}=\Gamma\Big(L\Delta Q_{\alpha\beta}-aQ_{\alpha\beta}+b[Q_{\alpha\gamma}Q_{\gamma\beta}-\frac{\delta_{\alpha\beta}}{d}\textrm{tr}(Q^2)]-cQ_{\alpha\beta}\textrm{tr}(Q^2)\Big) \\
      \partial_t u_\alpha+u_\beta\partial_\beta u_\alpha=\nu\Delta u_\alpha+\partial_\alpha p-L \partial_\beta \left(\partial_\alpha Q_{\zeta\delta}\partial_\beta Q_{\zeta\delta}-\frac{\delta_{\alpha\beta}}{d}\partial_\lambda Q_{\zeta\delta}\partial_\lambda Q_{\zeta\delta}\right)
      +L\partial_\beta\left(Q_{\alpha\gamma}\Delta Q_{\gamma\beta}-\Delta Q_{\alpha\gamma}Q_{\gamma\beta}\right) \\
       \partial_\gamma u_\gamma=0
     \end{array}\right.
      \label{system}
\end{equation}  in $\mathbb{R}^d$, $d=2,3$.

\par We also need to assume from now on that 
\begin{equation}
c>0
\label{c+}
\end{equation} This assumption is necessary from a modelling point of view (see \cite{apala},\cite{mz} ) so that the energy $\mathcal{F}$ (see next section, relation (\ref{freeenergy})) is bounded from below, and it is also necessary for having global solutions (see Proposition ~\ref{prop:aprioriest} and its proof).

\par We restrict ourselves to the case  $\xi=0$ for technical simplicity. However,  we think that our method can also be used in the general case $\xi\not=0$ and we will study this in a forthcoming paper \cite{pz2}.
\smallskip
\par The paper is organised as follows: in the second section we show that the equation admits a Lyapunov functional, whose existence is based on a certain cancellation that will prove to be crucial in the proof of higher regularity and the weak-strong uniqueness. Using the apriori estimates provided by the existence of a Lyapunov functional we show in the third section the existence of weak solutions in dimensions two and three. In the fourth section we restrict to dimension two and show the existence of arbitrarily regular solutions and that the strong norms  increase in time at most triply exponentially. Finally in the last section we show the weak-strong uniqueness of solutions in dimension two. The appendix contains a technical calculation necessary in the fourth section.

\smallskip\par{\bf Notations and conventions}  Let
$S_0\subset \mathbb{M}^{3\times 3}$ denote the space of Q-tensors,  i.e.
 $$S_0\stackrel{def}{=} \left\{Q \in \mathbb{M}^{3\times 3};
Q_{ij}=Q_{ji},\textrm{tr}(Q) = 0, i,j=1,2,3 \right\}$$
\par We use the Frobenius norm of a matrix $
 \left| Q \right|\stackrel{def}{=}\sqrt{\textrm{tr}Q^2} =\sqrt{ Q_{\alpha\beta}
Q_{\alpha\beta}}$ and define Sobolev spaces of $Q$-tensors in terms of this norm. For instance $H^1(\mathbb{R}^d,S_0)\stackrel{def}{=}\{Q:\mathbb{R}^d\to S_0, \int_{\mathbb{R}^d} |\nabla Q(x)|^2+|Q(x)|^2\,dx<\infty\}$. For $A,B\in S_0$ we denote $A\cdot B=\textrm{tr}(AB)$ and $|A|=\sqrt{\textrm{tr}(A^2)}$. We also denote  $|\nabla Q|^2(x)\stackrel{def}{=}Q_{\alpha\beta,\gamma}(x)Q_{\alpha\beta,\gamma}(x)$ and $|\Delta Q|^2(x)\stackrel{def}{=}\Delta Q_{\alpha\beta}(x)\Delta Q_{\alpha\beta}(x)$. We recall also that $\Omega_{\alpha\beta}\stackrel{def}{=}\frac{1}{2}\left(\partial_\beta u_\alpha-\partial_\alpha u_\beta\right)$ and $u_{\alpha,\beta}\stackrel{def}{=}\partial_\beta u_\alpha$, $Q_{ij,k}\stackrel{def}{=}\partial_k Q_{ij}$.

\section{ The dissipation principle and apriori estimates}

\par Let us denote the free energy of the director fields:

\begin{equation}
\mathcal{F}(Q)=\int_{\mathbb{R}^d} \frac{L}{2}|\nabla Q|^2+\frac{a}{2}\textrm{tr}(Q^2)-\frac{b}{3}\textrm{tr}(Q^3)+\frac{c}{4}\textrm{tr}^2(Q^2)\,dx
\label{freeenergy}
\end{equation}

\par In the absence of the flow, when $u=0$ in the equations (\ref{system}), the free energy is a Lyapunov functional of the system. If $u\not=0$ we still have a Lyapunov functional for (\ref{system}) but this time one that includes the kinetic energy of the system. More precisely we have:

\begin{proposition}
The system (\ref{system}) has a Lyapunov functional:
\begin{equation}
E(t)\stackrel{def}{=}\frac{1}{2}\int_{\mathbb{R}^d}|u|^2(t,x)\,dx+\int_{\mathbb{R}^d}\frac{L}{2} |\nabla Q|^2(t,x)+\frac{a}{2}\textrm{tr}(Q^2(t,x))-\frac{b}{3}\textrm{tr}(Q^3(t,x))+\frac{c}{4}\textrm{tr}^2(Q^2(t,x))\,dx
\end{equation}
\par If $d=2,3$ and $(Q,u)$ is a smooth solution of (\ref{system}) such that $Q\in L^\infty(0,T; H^1(\mathbb{R}^d))\cap L^2(0,T;H^2(\mathbb{R}^d))$  and $u\in L^\infty(0,T;L^2(\mathbb{R}^d))\cap L^2(0,T;H^1(\mathbb{R}^d))$ then, for all $t<T$, we have:
\begin{equation}
\frac{d}{dt}E(t)=-\nu\int_{\mathbb{R}^d}|\nabla u|^2\,dx-\Gamma\int_{\mathbb{R}^d} \textrm{tr}\left(L\Delta Q-aQ+b[Q^2-\frac{\textrm{tr}(Q^2)}{3}Id]-cQ\textrm{tr}(Q^2)\right)^2\,dx\le 0
\label{energydecay}
\end{equation}
\label{prop:Lyapunov}
\end{proposition}

\smallskip\par{\bf Proof.} We multiply the first equation in (\ref{system}) to the right by $-\left(L\Delta Q-aQ+b[Q^2-\frac{\textrm{tr}(Q^2)}{3}Id]-cQ\textrm{tr}(Q^2)\right)$, take the trace, integrate over $\mathbb{R}^d$ and by parts and sum with the second equation multiplied by $u$ and integrated over $\mathbb{R}^d$ and by parts (let us observe that because of our assumptions on $Q$ and $u$ we do not have boundary terms, when integrating by parts). We obtain:

\begin{eqnarray}
\frac{d}{dt}\int_{\mathbb{R}^d}\frac{1}{2}|u|^2+\frac{L}{2}|\nabla Q|^2+\frac{a}{2}\textrm{tr}(Q^2)-\frac{b}{3}\textrm{tr}(Q^3)+\frac{c}{4}\textrm{tr}^2(Q^2)\,dx\nonumber\\+\nu\int_{\mathbb{R}^d}|\nabla u|^2\,dx+\Gamma\int_{\mathbb{R}^d}\textrm{tr}\left(L\Delta Q-aQ+b[Q^2-\frac{\textrm{tr}(Q^2)}{3}Id]-cQ\textrm{tr}(Q^2)\right)^2\,dx\nonumber\\=\underbrace{\int_{\mathbb{R}^d}u\cdot\nabla Q_{\alpha\beta}\left(-aQ_{\alpha\beta}+b[Q_{\alpha\gamma}Q_{\gamma\beta}-\frac{\delta_{\alpha\beta}}{3}\textrm{tr}(Q^2)]-cQ_{\alpha\beta}\textrm{tr}(Q^2))\right)\,dx}_{\mathcal{I}}\nonumber\\+\underbrace{\int_{\mathbb{R}^d}\left(-\Omega_{\alpha\gamma} Q_{\gamma\beta}+Q_{\alpha\gamma}\Omega_{\gamma\beta}\right)\left(-aQ_{\alpha\beta}+b[Q_{\alpha\delta}Q_{\delta\beta}-\frac{\delta_{\alpha\beta}}{3}\textrm{tr}(Q^2)]-cQ_{\alpha\beta}\textrm{tr}(Q^2))\right)\,dx}_{\mathcal{II}}
\nonumber\\
+\underbrace{L\int_{\mathbb{R}^d}u_{\gamma} Q_{\alpha\beta,\gamma}\Delta Q_{\alpha\beta}\,dx}_{\mathcal{A}}\underbrace{
-\frac{L}{2}\int_{\mathbb{R}^d} u_{\alpha,\gamma}Q_{\gamma\beta}\Delta Q_{\alpha\beta}\,dx}_{\mathcal{B}}
\nonumber
\end{eqnarray}
\begin{eqnarray}
\underbrace{+\frac{L}{2}\int_{\mathbb{R}^d}u_{\gamma,\alpha}Q_{\gamma\beta}\Delta Q_{\alpha\beta}}_{\mathcal{C}}\,dx
\underbrace{+\frac{L}{2}\int_{\mathbb{R}^d}Q_{\alpha\gamma}u_{\gamma,\beta}\Delta Q_{\alpha\beta}\,dx}_{\mathcal{C}}\underbrace{-\frac{L}{2}\int_{\mathbb{R}^d}Q_{\alpha\gamma}u_{\beta,\gamma}\Delta Q_{\alpha\beta}\,dx}_{\mathcal{B}}\nonumber\\
\underbrace{+L\int_{\mathbb{R}^d}Q_{\gamma\delta,\alpha}Q_{\gamma\delta,\beta}u_{\alpha,\beta}\,dx}_{\mathcal{AA}}-L\int_{\mathbb{R}^d} Q_{\alpha\gamma}\Delta Q_{\gamma\beta}u_{\alpha,\beta}\,dx+L\int_{\mathbb{R}^d}\Delta Q_{\alpha\gamma}Q_{\gamma\beta}u_{\alpha,\beta}\,dx\nonumber
\end{eqnarray}
\begin{eqnarray}
=\underbrace{-L\int_{\mathbb{R}^d} u_{\alpha,\gamma}Q_{\gamma\beta}\Delta Q_{\alpha\beta}\,dx}_{2\mathcal{B}}\underbrace{+L\int_{\mathbb{R}^d}u_{\gamma,\alpha}Q_{\gamma\beta}\Delta Q_{\alpha\beta}\,dx}_{2\mathcal{C}}
\nonumber\\ \underbrace{-L\int_{\mathbb{R}^d} Q_{\alpha\gamma}\Delta Q_{\gamma\beta}u_{\alpha,\beta}\,dx}_{\mathcal{CC}}\underbrace{+L\int_{\mathbb{R}^d}\Delta Q_{\alpha\gamma}Q_{\gamma\beta}u_{\alpha,\beta}\,dx}_{\mathcal{BB}}=0
\label{Lyapunovcancellation}
\end{eqnarray} where $\mathcal{I}=0$ (since $\nabla\cdot u=0$), $\mathcal{II}=0$ (since $Q_{\alpha\beta}=Q_{\beta\alpha}$) and for the second equality we used

\begin{eqnarray}\underbrace{\int_{\mathbb{R}^d}u_\gamma Q_{\alpha\beta,\gamma}\Delta Q_{\alpha\beta}\,dx}_{\mathcal{A}}\underbrace{+\int_{\mathbb{R}^d}Q_{\gamma\delta,\alpha}Q_{\gamma\delta,\beta}u_{\alpha,\beta}\,dx}_{\mathcal{AA}}=\int_{\mathbb{R}^d}u_\gamma Q_{\alpha\beta,\gamma}\Delta Q_{\alpha\beta}\,dx\nonumber\\-\int_{\mathbb{R}^d}Q_{\gamma\delta,\alpha}Q_{\gamma\delta,\beta\beta}u_\alpha\,dx-
\int_{\mathbb{R}^d}Q_{\gamma\delta,\alpha\beta}Q_{\gamma\delta,\beta}u_\alpha\,dx=\int_{\mathbb{R}^d}\frac{1}{2}Q_{\gamma\delta,\alpha}Q_{\gamma\delta,\alpha}u_{\alpha,\alpha}\,dx=0
\end{eqnarray} while for the last equality in (\ref{Lyapunovcancellation}) we used $2\mathcal{B}+\mathcal{BB}=2\mathcal{C}+\mathcal{CC}=0$. $\Box$

\bigskip
\par In the following we assume that there exists a smooth solution of (\ref{system}) and  obtain estimates on the behaviour of various norms:

\begin{proposition} Let $(Q,u)$ be a smooth solution of (\ref{system}), with restriction (\ref{c+}), and smooth initial data $(\bar Q(x),\bar u(x))$, that decays fast enough at infinity so that we can integrate by parts in space (for any $t\ge 0$) without boundary terms.
\par (i) If $\bar Q\in L^p$ for some $p\ge 2$ we have
\begin{equation}
\|Q(t,\cdot)\|_{L^p}\le e^{Ct}\|\bar Q\|_{L^p}, \forall t\ge 0
\end{equation} with $C=C(a,b,c,p,\Gamma)$.
\par (ii) For $d=2,3$  (and $(\bar Q,\bar u)$ so that the right hand side of the expression below is finite) we have:
\begin{eqnarray}
\|u(t,\cdot)\|_{L^2}^2+2\nu\int_0^t\|\nabla u(s,\cdot)\|_{L^2}^2\,ds+L\|\nabla Q(t,\cdot)\|_{L^2}^2+\Gamma L^2\int_0^t \|\Delta Q(s,\cdot)\|_{L^2}^2\,ds\le \|u(0,\cdot)\|_{L^2}^2+\|\nabla Q(0,\cdot)\|_{L^2}^2\nonumber\\+Ce^{Ct}\left(\|Q(0,\cdot)\|_{L^2}^2+\|Q(0,\cdot)\|_{L^6}^6\right)
\end{eqnarray} with the constant $C=C(a,b,c,d,L,\Gamma)$.
\label{prop:aprioriest}
\end{proposition}

\smallskip\par{\bf Proof.}
\par (i) Multiplying the first equation in (\ref{system}) by $2pQ\textrm{tr}^{p-1}(Q^2)$ and taking the trace we obtain:

\begin{eqnarray}
\left(\partial_t+u\cdot\nabla\right)\textrm{tr}^p(Q^2)=\Gamma\Big(2pL\Delta Q_{\alpha\beta}Q_{\alpha\beta}\textrm{tr}^{p-1}(Q^2)-2pa\textrm{tr}^p(Q^2)+2pb\textrm{tr}(Q^3)\textrm{tr}^{p-1}(Q^2)-2pc\textrm{tr}^{p+1}(Q^2)\Big)
\end{eqnarray}

\par Let us observe that for $Q$  a traceless, symmetric, $3\times 3$ matrix we have:

\begin{equation}
\textrm{tr}(Q^3)\le \frac{3\varepsilon}{8}\textrm{tr}^2(Q^2)+\frac{1}{\varepsilon} \textrm{tr}(Q^2),\forall \varepsilon>0
\label{est:Q3}
\end{equation}

\par Indeed, if $Q$ has the eigenvalues $x,y,-x-y$ then $\textrm{tr}(Q^3)=-3xy(x+y)$, $\textrm{tr}(Q^2)=2(x^2+y^2+xy)$ and the inequality (\ref{est:Q3}) follows.

\par Integrating over $\mathbb{R}^d$, integrating by parts ( we have no boundary terms because of our assumption), as well as using that $\nabla\cdot u=0$, together with (\ref{est:Q3})(where $\varepsilon=\frac{4c}{3b}$) and  the assumption $c>0$ we obtain:

\begin{eqnarray}
\partial_t\int_{\mathbb{R}^d}\textrm{tr}^p(Q^2)\,dx\le \underbrace{-2p \Gamma L\int_{\mathbb{R}^d}\nabla Q_{\alpha\beta}\nabla Q_{\alpha\beta}\textrm{tr}^{p-1}(Q^2))\,dx}_{\le 0}\nonumber\\
\underbrace{-4p(p-1)\Gamma L\int_{\mathbb{R}^d}Q_{\alpha\beta,\gamma}Q_{\alpha\beta}Q_{\delta\lambda,\gamma}Q_{\delta\lambda}\textrm{tr}^{p-2}(Q^2)\,dx}_{\le 0}
+C\int_{\mathbb{R}^d}\textrm{tr}^p(Q^2)\,dx
\end{eqnarray} where the constant $C$ depends on $a,b,c, p$ and $\Gamma$. Thus we have

\begin{equation}
\int_{\mathbb{R}^d}\textrm{tr}^p(Q^2(t,x))\,dx\le e^{Ct}\int_{\mathbb{R}^d}\textrm{tr}^p(Q^2(0,x))\,dx
\label{Lpest}
\end{equation} with $C=C(a,b,c,p,\Gamma)$.

\smallskip\par (ii) Relation (\ref{energydecay}) implies

\begin{eqnarray}
 \frac{L}{2}\|\nabla Q(t,\cdot)\|_{L^2}^2+\frac{1}{2}\|u(t,\cdot)\|_{L^2}^2+\nu\int_0^t \|\nabla u(s,\cdot)\|_{L^2}^2\,ds+\Gamma L^2 \int_0^t \|\Delta Q(s,\cdot)\|_{L^2}^2ds\nonumber\\
 \le C\int_{\mathbb{R}^d} \textrm{tr}(Q^2(t,x))+\textrm{tr}^2(Q^2(t,x))\,dx+C\int_{\mathbb{R}^d} \textrm{tr}(Q^2(0,x))+\textrm{tr}^2(Q^2(0,x))\,dx+ \frac{L}{2}\|\nabla Q(0,\cdot)\|_{L^2}+\frac{1}{2}\|u(0,\cdot)\|_{L^2}^2\nonumber\\
 + \Gamma \int_0^t\int_{\mathbb{R}^d} \textrm{tr}\Big(L\Delta Q\big(aQ-bQ^2+cQ\textrm{tr}(Q^2)\big)\Big)\,dx\,ds+\Gamma \int_0^t\int_{\mathbb{R}^d} \textrm{tr}\Big(\big(aQ-bQ^2+cQ\textrm{tr}(Q^2)\big)L\Delta Q\Big)\,dx\,ds\nonumber
 \end{eqnarray}

\par In the last inequality we  use Holder inequality to estimate $\Delta Q$ in $L^2$  and absorb it in the left hand side  while the terms without gradients are estimated using (\ref{Lpest}) and interpolation between the $L^2$ and $L^6$ norms.$\Box$

\section{Weak solutions}

\par A pair $(Q,u)$ is called a weak solution  of the system (\ref{system}), subject to initial data
\begin{equation}
Q(0,x)=\bar Q(x)\in L^2(\mathbb{R}^d),\, u(0,x)=\bar u(x)\in L^2(\mathbb{R}^d), \nabla\cdot \bar u=0\,\textrm{ in }\mathcal{D}'(\mathbb{R}^d)
\label{initialdata}
\end{equation}
if   $Q\in L^\infty_{loc}(\mathbb{R}_+;H^1)\cap L^2_{loc}(\mathbb{R}_+;H^2)$, $u\in L^\infty_{loc}(\mathbb{R}_+;L^2)\cap L^2_{loc}(\mathbb{R}_+;H^1)$  and for every compactly supported $\varphi\in C^\infty([0,\infty)\times \mathbb{R}^d; S_0)$, $\psi\in C^\infty ([0,\infty)\times\mathbb{R}^d;\mathbb{R}^d)$ with $\nabla\cdot\psi=0$ we have
\begin{eqnarray}\int_0^\infty\int_{\mathbb{R}^d}(-Q\cdot \partial_t\varphi-\Gamma L\Delta Q\cdot \varphi)- Q\cdot u\nabla_x \varphi-\Omega
Q\cdot\varphi +Q\Omega\cdot\varphi \,\,dx\,dt\nonumber\\
=\int_{\mathbb{R}^d}\bar Q(x)\cdot \varphi(0,x)\,dx+\Gamma\int_0^\infty\int_{\mathbb{R}^d}\Big\{-aQ+b[Q^2-\frac{\textrm{tr}(Q^2)}{d}Id]-cQ\textrm{tr}(Q^2)\Big\}\cdot\varphi\,\,dx\,dt
\label{weaksol1}
\end{eqnarray} 
and
\begin{eqnarray}\int_0^\infty\int_{\mathbb{R}^d}-u\partial_t\psi-u_\alpha u_\beta\partial_\alpha\psi_\beta+\nu\nabla u\nabla\psi\,\, dt\,dx-\int_{\mathbb{R}^d}\bar u(x)\psi(0,x)\,dx\nonumber\\
=L\int_0^\infty\int_{\mathbb{R}^d} 
 Q_{\gamma\delta,\alpha} Q_{\gamma\delta,\beta}\psi_{\alpha,\beta} -Q_{\alpha\gamma}\Delta Q_{\gamma\beta}\psi_{\alpha,\beta}+\Delta Q_{\alpha\gamma}Q_{\gamma\beta}\psi_{\alpha,\beta}\,\,dx\,dt
\label{weaksol2}
\end{eqnarray} 
\begin{proposition} For $d=2,3$ there exists a weak solution $(Q,u)$ of the system (\ref{system}), with restriction (\ref{c+}), subject to initial conditions (\ref{initialdata}). The solution $(Q,u)$ is such that $Q\in L^\infty_{loc}(\mathbb{R}_+;H^1)\cap L^2_{loc}(\mathbb{R}_+;H^2)$ and $u\in L^\infty_{loc}(\mathbb{R}_+;L^2)\cap L^2_{loc}(\mathbb{R}_+;H^1)$.
\label{prop:weak}
\end{proposition} 

\smallskip\par{\bf Proof.}  We define the mollifying operator

$$\widehat{J_nf}(\xi)=1_{[\frac{1}{n},n]}(|\xi|)\hat f(\xi)$$ and consider the system:

\begin{equation}
\left\{\begin{array}{l}
      \partial_t Q^{(n)}+J_n \Big(\mathcal{P}J_n u^n \nabla J_nQ^{(n)}\Big)-J_n\Big(\mathcal{P}J_n\Omega^n J_n Q^{(n)}\Big)+J_n\Big(J_n Q^{(n)}\mathcal{P}J_n\Omega^n\Big)=\Gamma L\Delta J_n Q^{(n)}\\
      +\Gamma\Big(-aJ_n Q^{(n)}+b[J_n(J_nQ^{(n)} J_n Q^{(n)})-\frac{\textrm{tr}(J_n(J_nQ^{(n)}J_n Q^{(n)}))}{d}Id]-cJ_nQ^{(n)}\textrm{tr}(J_n(J_nQ^{(n)}J_n Q^{(n)}))\Big) \\
      \partial_t u^n +\mathcal{P}J_n(\mathcal{P}J_nu^n\nabla\mathcal{P}J_n u^n)=- L\mathcal{P}J_n(\nabla\cdot\left(\textrm{tr}(\nabla J_n Q^{(n)}\nabla J_n Q^{(n)})-\frac{1}{d}|\nabla J_n Q^{(n)}|^2 Id\right))\nonumber\\+L\mathcal{P}(\nabla\cdot J_n\left(J_nQ^{(n)} \Delta J_nQ^{(n)}-\Delta J_n Q^{(n)} J_n Q^{(n)}\right))+\nu\Delta \mathcal{P}J_nu^n\\
     \end{array}\right.
      \label{approxsystem+}
\end{equation} where $\mathcal{P}$ denotes the Leray projector onto divergence-free vector fields.

\par The system above can be regarded as an ordinary differential equation in $L^2$ verifying the conditions of the Cauchy-Lipschitz theorem. Thus it admits a unique maximal solution $(Q^{(n)}, u^n)\in C^1([0,T_n); L^2(\mathbb{R}^d;\mathbb{R}^{d\times d})\times L^2(\mathbb{R}^d,\mathbb{R}^d))$. As we have $(\mathcal{P}J_n)^2=\mathcal{P}J_n$ and $J_n^2=J_n$ the pair $(J_n Q^{(n)}, \mathcal{P}J_n u^n)$ is also a solution of (\ref{approxsystem+}). By uniqueness we have $(J_n Q^{(n)}, \mathcal{P}J_n u^n)=(Q^{(n)}, u^n)$ hence $(Q^{(n)}, u^n)\in C^1([0,T_n),H^\infty)$ and $(Q^{(n)}, u^n)$ satisfy the system:

\begin{equation}
\left\{\begin{array}{l}
      \partial_t Q^{(n)}+J_n\big(u^n \nabla Q^{(n)}\big)-J_n\big(\Omega^n  Q^{(n)}-Q^{(n)}\Omega^n\big)=\Gamma L\Delta Q^{(n)}\\
      +\Gamma\Big(-a Q^{(n)}+b[J_n(Q^{(n)}  Q^{(n)})-\frac{\textrm{tr}(J_n(Q^{(n)} Q^{(n)}))}{d}Id]-cQ^{(n)}\textrm{tr}(J_n(Q^{(n)} Q^{(n)}))\Big) \\
      \partial_t u^n +\mathcal{P}J_n(u^n\nabla u^n)=- L\mathcal{P}J_n(\nabla\cdot\left(\textrm{tr}(\nabla  Q^{(n)}\nabla  Q^{(n)})-\frac{1}{d}|\nabla Q^{(n)}|^2 Id\right))\\+L\mathcal{P}(\nabla\cdot J_n\left(Q^{(n)} \Delta Q^{(n)}-\Delta Q^{(n)} Q^{(n)} \right))+\nu\Delta u^n
     \end{array}\right.
      \label{approxsystem++}
\end{equation} 

\par We can argue as in the proof of  the apriori estimates and the same estimates hold for the approximating system (\ref{approxsystem++}). These estimates allow us to conclude that $T_n=\infty$ and we also get the following apriori bounds:

\begin{eqnarray}
\sup_n \|Q^{(n)}\|_{L^2(0,T;H^2)\cap L^\infty(0,T;H^1)}<\infty\nonumber\\
\sup_n \|u^n\|_{L^\infty(0,T;L^2)\cap L^2(0,T;H^1)}<\infty
\label{weaksolapriori}
\end{eqnarray} for any $T<\infty$.

\par The pair $(Q^{(n)}, u^n)$ is also a weak solution of the approximating system (\ref{approxsystem++}) hence for every compactly supported $\varphi\in C^\infty([0,\infty)\times \mathbb{R}^d; S_0)$, $\psi\in C^\infty ([0,\infty)+\times\mathbb{R}^d;\mathbb{R}^d)$ with $\nabla\cdot\psi=0$ we have:

\begin{eqnarray}\int_0^\infty\int_{\mathbb{R}^d}(-Q^{(n)}\cdot \partial_t\varphi-\Gamma L \Delta Q^{(n)}\cdot\varphi) - J_n\big(Q^{(n)}\cdot u^n\big)\nabla_x\varphi-J_n\big(\Omega^n
Q^{(n)}\big)\cdot\varphi +J_n\big(Q^{(n)}\Omega^n\big)\cdot\varphi \,\,dx\,dt\nonumber\\
=\int_{\mathbb{R}^d}\bar Q(x)\cdot \varphi(0,x)\,dx+\Gamma\int_0^\infty\int_{\mathbb{R}^d}\{-aQ^{(n)}+b[J_n\big(\left(Q^{(n)}\right)^2\big)-\frac{\textrm{tr}\big(J_n(\left(Q^{(n)}\right)^2)\big)}{d}Id]-cQ^{(n)}\textrm{tr}(J_n(Q^{(n)})^2)\}\cdot\varphi\,\,dx\,dt
\label{weaksol1+}
\end{eqnarray}
and
\begin{eqnarray}\int_0^\infty\int_{\mathbb{R}^d}-u^n\partial_t\psi-J_n(u^n_\alpha u^n_\beta)\partial_\alpha\psi_\beta+\nu\nabla u^n\nabla\psi\,\,dx\, dt-\int_{\mathbb{R}^d}\bar u(x)\psi(0,x)\,dx\nonumber\\
=L\int_0^\infty\int_{\mathbb{R}^d}\Big\{ J_n\left(Q^{(n)}_{\gamma\delta,\alpha} Q^{(n)}_{\gamma\delta,\beta}\right)\psi_{\alpha,\beta}-J_n\left(Q^{(n)}_{\alpha\gamma}\Delta Q^{(n)}_{\gamma\beta}-\nu\Delta Q^{(n)}_{\alpha\gamma}Q^{(n)}_{\gamma\beta}\right)\psi_{\alpha,\beta}\Big\}\,\,dx\,dt
\label{weaksol2+}
\end{eqnarray}

\par We consider the solutions of (\ref{approxsystem++}) and taking into account the bounds (\ref{weaksolapriori}) we get, by classical compactness and weak convergence arguments, that there exists a $Q\in L^\infty_{loc}(\mathbb{R}_+;H^1)\cap L^2_{loc}(\mathbb{R}_+;H^2)$ and a $u\in L^\infty_{loc}(\mathbb{R}_+;L^2)\cap L^2_{loc}(\mathbb{R}_+;H^1)$ so that, on a subsequence, we have:

\begin{eqnarray}
Q^{(n)}\rightharpoonup Q\textrm{ in } L^2(0,T;H^2)\,\textrm{ and }Q^{(n)}\to Q\textrm{ in }L^2(0,T;H_{loc}^{2-\varepsilon}),\forall \varepsilon>0\nonumber\\
Q^{(n)}(t)\rightharpoonup Q(t)\textrm{ in }H^1\textrm{ for all }t\in\mathbb{R}_+\nonumber\\
u^n\rightharpoonup u\textrm{ in }L^2(0,T;H^1)\,\textrm{ and }u^n\to u\textrm{ in }L^2(0,T;H_{loc}^{1-\varepsilon}),\forall\varepsilon>0\nonumber\\
u^n(t)\rightharpoonup u(t)\,\textrm{ in }L^2\textrm{ for all }t\in\mathbb{R}_+
\label{convergences}
\end{eqnarray}

 These  convergences allow us to the pass to the limit in the weak solutions 
(\ref{weaksol1+}),(\ref{weaksol2+}) to obtain a weak solution of (\ref{system}), namely  (\ref{weaksol1}),(\ref{weaksol2}). The term that is the most difficult to treat in passing  to the limit is  the last term in (\ref{weaksol2+}), namely 

$$L\int_0^\infty\int_{\mathbb{R}^d}J_n\left(Q^{(n)}_{\alpha\gamma}\Delta Q^{(n)}_{\gamma\beta}-\Delta Q^{(n)}_{\alpha\gamma}Q^{(n)}_{\gamma\beta}\right)\psi_{\alpha,\beta}\,\,dx\,dt=L\int_0^\infty\int_{\mathbb{R}^d}\left(Q^{(n)}_{\alpha\gamma}\Delta Q^{(n)}_{\gamma\beta}-\Delta Q^{(n)}_{\alpha\gamma}Q^{(n)}_{\gamma\beta}\right)\cdot J_n\psi_{\alpha,\beta}\,\,dx\,dt.$$

 Recalling that $\psi$ is compactly supported we have that there exists a time $T>0$ so that $\psi(t,x)=J_n\psi(t,x)=0,\forall t>T, x\in\mathbb{R}^d,n\in\mathbb{N}$. Taking into account that $\psi$ is compactly supported and the convergences (\ref{convergences}) one can easily pass to the limit the terms $\partial_\beta J_n\psi_\alpha Q^{(n)}_{\alpha\gamma}$ and $\partial_\beta J_n\psi_\alpha Q^{(n)}_{\gamma\beta}$ strongly in $L^2(0,T;L^2)$. Indeed we have:

\begin{equation}
\partial_\beta J_n\psi_\alpha Q^{(n)}_{\alpha\gamma}-\partial_\beta\psi_\alpha Q_{\alpha\gamma}=
\underbrace{\Big(\partial_\beta J_n\psi_\alpha-\partial_\beta\psi_\alpha\Big) Q^{(n)}_{\alpha\gamma}}_{\mathcal{I}}+
\underbrace{\partial_\beta\psi_\alpha\Big(Q^{(n)}_{\alpha\gamma}-Q_{\alpha\gamma}\Big)}_{\mathcal{II}}
\end{equation} and the first term, $\mathcal{I}$, converges to $0$, strongly in $L^2(0,T;L^2)$ because $\psi$ is smooth and compactly supported, hence $\partial_\beta J_n\psi-\partial_\beta\psi$ converges to zero in any $L^q(0,T;L^p)$  and $Q^{(n)}$ is bounded in $L^\infty$ in time and $L^p$ in space ($1<p<\infty$ if $d=2$ and $2\le p\le 6$ if $d=3$, due to the bounds (\ref{weaksolapriori})). On the other hand  the second term $\mathcal{II}$ converges strongly to zero in $L^2(0,T;L^2)$ because of (\ref{convergences}) and the fact that $\psi$ is compactly supported.
 
 \par Relations (\ref{convergences}) give that   $\Delta Q^{(n)}_{\gamma\beta}$, $\Delta Q^{(n)}_{\alpha\gamma}$ converges weakly in $L^2(0,T;L^2)$. Thus  we get convergence to the limit term 
 
 \begin{eqnarray}
 L\int_0^\infty \int_{\mathbb{R}^d} (\Delta Q_{\gamma\beta})(\partial_\beta \psi_\alpha Q_{\alpha\gamma})dxdt-L\int_0^\infty\int_{\mathbb{R}^d}(\Delta Q_{\alpha\gamma})(\partial_\beta\psi_\alpha Q_{\gamma\beta})dxdt\nonumber\\
=L\int_0^T \int_{\mathbb{R}^d} (\Delta Q_{\gamma\beta})(\partial_\beta \psi_\alpha Q_{\alpha\gamma})dxdt-L\int_0^T\int_{\mathbb{R}^d}(\Delta Q_{\alpha\gamma})(\partial_\beta\psi_\alpha Q_{\gamma\beta})dxdt.
 \end{eqnarray} $\Box$

\section{ Higher regularity in $2$D, using the dissipation principle}

\par In this section we restrict ourselves to dimension two and show that starting from an initial data with some higher regularity, we can obtain more regular solutions. More precisely, we have:

\begin{theorem}
Let $s>1$ and $(\bar Q,\bar u)\in H^{s+1}(\mathbb{R}^2)\times H^s(\mathbb{R}^2)$. There exists a global a solution $(Q(t,x),u(t,x))$ of the system (\ref{system}), with restriction (\ref{c+}), subject  to initial conditions
$$Q(0,x)=\bar Q(x),\, u(0,x)=\bar u(x)$$ and $Q\in L^2_{loc}(\mathbb{R}_+; H^{s+2}(\mathbb{R}^2))\cap L^\infty_{loc}(\mathbb{R}_+;H^{s+1}(\mathbb{R}^2))$, $u\in L^2_{loc}(\mathbb{R}_+;H^{s+1}(\mathbb{R}^2)\cap L^\infty_{loc}(\mathbb{R}_+;H^s)$.
\par Moreover, we have:

\begin{equation}
L\|\nabla Q(t,\cdot)\|_{H^s(\mathbb{R}^2)}^2+\|u(t,\cdot)\|_{H^s(\mathbb{R}^2)}^2\le \Big(e+\|\bar Q\|_{H^{s+1}(\mathbb{R}^2)}+\|\bar u\|_{H^s(\mathbb{R}^2)}\Big)^{e^{e^{Ct}}}
\label{rate:3exp}
\end{equation}  where the constant $C$ depends only on $\bar Q, \bar u$, $a,b,c$, $\Gamma$ and $L$. 
\label{theorem:reg}
\end{theorem}

\par The proof of the theorem is mainly based on $H^s$ energy estimates and the following cancelation(that is also used implicitly in showing the dissipation of the energy in Proposition ~\ref{prop:Lyapunov}): 

\begin{lemma}
\label{anulare}
For any symmetric matrices $Q', Q\in\mathbb{R}^{d\times d}$ and $\Omega_{\alpha\beta}=\frac{1}{2}(u_{\alpha,\beta}-u_{\beta,\alpha})\in\mathbb{R}^{d\times d}$ we have
 $$\int_{\mathbb{R}^d} \textrm{tr}\big((\Omega Q' -Q'\Omega)\Delta Q\big)\,dx-\int_{\mathbb{R}^d}\partial_\beta(Q'_{\alpha\gamma}\Delta
Q_{\gamma\beta}-\Delta Q_{\alpha\gamma} Q'_{\gamma\beta})u_\alpha\,dx=0$$
\end{lemma}

\smallskip\par {\bf Proof.} We note that

\begin{eqnarray}\int_{\mathbb{R}^d} \textrm{tr}\big((\Omega Q' -Q'\Omega)\Delta Q\big)\,dx=\int_{\mathbb{R}^d}\Omega_{\alpha\gamma}Q'_{\gamma\beta}\Delta Q_{\beta\alpha}-Q'_{\alpha\gamma}\Omega_{\gamma\beta}\Delta Q_{\beta\alpha}
=\int_{\mathbb{R}^d}\Omega_{\alpha\gamma}Q'_{\gamma\beta}\Delta Q_{\beta\alpha}+\Omega_{\beta\gamma}
Q'_{\gamma\alpha}\Delta Q_{\alpha\beta}\nonumber\\=2\int_{\mathbb{R}^d}\textrm{tr}\big(\Omega Q'\Delta Q\big)\,dx=\underbrace{\int_{\mathbb{R}^d} u_{\alpha,\beta}Q'_{\beta\gamma}\Delta Q_{\gamma\alpha}\,dx}_{\mathcal{I}_1}-\underbrace{\int_{\mathbb{R}^d}
u_{\beta,\alpha}Q'_{\beta\gamma}\Delta Q_{\gamma\alpha}\,dx}_{\mathcal{I}_2}
\end{eqnarray} and on the other hand 
$$-\int_{\mathbb{R}^d} \partial_\beta (Q'_{\alpha\gamma}\Delta Q_{\gamma\beta}) u_\alpha=\int_{\mathbb{R}^d} Q'_{\alpha\gamma}\Delta Q_{\gamma \beta}\partial_\beta u_\alpha=\int_{\mathbb{R}^d} Q'_{\beta\gamma}\Delta Q_{\gamma \alpha}\partial_\alpha u_\beta=I_2$$ 
and also 
$$\int_{\mathbb{R}^d} \partial_\beta(\Delta Q_{\alpha\gamma} Q'_{\gamma\beta}) u_\alpha=-\int_{\mathbb{R}^d}  Q'_{\beta \gamma}\Delta Q_{\gamma \alpha}\partial_\beta u_\alpha=-I_1$$
which finishes the proof. $\Box$

\smallskip\par\begin{remark} The main point in the proof of the theorem is to use the previous lemma to eliminate the highest derivatives in $u$ in the first equation of the system (\ref{system}) and the highest derivatives in $Q$ in the second equation of the system.The proof  could have been done, alternatively, by differentiating the equations $k\ge 1$ times and using the previous lemma. However that would have required estimating some delicate commutators and  would have restricted the initial data to $(\bar Q,\bar u)\in H^3\times H^2$. The Littlewood-Paley approach that we use allows for $(\bar Q,\bar u)\in H^{s+1}\times H^s$ with $s>1$.
\end{remark}

\smallskip
\par In order to prove the theorem we need to introduce some technical preliminaries:

\subsection{Littlewood-Paley theory}
We  define
$\mathcal{C}$ to be the ring of center
$0$, of small radius $1/2$ and great radius $2$. There exist two
nonnegative  radial
functions $\chi$ and $\varphi$ belonging respectively to~${\mathcal{D}} 
(B(0,1)) $ and to
${\mathcal{D}} (\mathcal{C}) $ so that
\begin{equation}
\label{lpfond1}
\chi(\xi) + \sum_{q\geq 0} \varphi (2^{-q}\xi) = 1,\forall \xi\in\mathbb{R}^d
\end{equation}
\begin{equation}
\label{lpfond2}
|p-q|\geq 2
\Rightarrow
{\rm Supp}\,\, \varphi(2^{-q}\cdot)\cap {\rm Supp}\,\, \varphi(2^{-p}\cdot)=\emptyset.
\end{equation}

\noindent
For instance, one can take $\chi \in \mathcal{D} (B(0,1))$ such that $
\chi  \equiv 1 $ on $B(0,1/2)$ and take
$$
\varphi(\xi) = \chi(\xi/2) -
\chi(\xi).
$$
Then, we are able to define the Littlewood-Paley decomposition. Let us denote
by~$\mathcal{F}$ the Fourier transform on~$\mathbb{R}^d$. Let
$h,\
\tilde h,\  \Delta_q, S_q$ ($q \in \mathbb{Z}$) be defined as follows:

$$\displaylines{
\label{defnotationdyadique}h = {\mathcal F}^{-1}\varphi\quad {\rm and}\quad \tilde h =
{\mathcal{F}}^{-1}\chi, \cr
\Delta_q u = \mathcal{F}^{-1}(\varphi(2^{-q}\xi)\mathcal{F} u) = 2^{qd}\int h(2^qy)u(x-y)dy,\cr
S_qu
=\mathcal{F}^{-1}(\chi(2^{-q}\xi)\mathcal{F} u) =2^{qd} \int \tilde h(2^qy)u(x-y)dy.\cr
}
$$

\par We recall that for two appropriately smooth functions $a$ and $b$ we have Bony's paraproduct decomposition
\cite{Bony81}:

\begin{equation}
ab=T_a b+T_b a+R(a,b)
\end{equation} where 
\begin{displaymath}
T_a b=\sum_{\substack{q'}}S_{q'-1} a\Delta_{q'}b,\,\, T_b a=\sum_{q'}S_{q'-1}b\Delta_{q'}a\textrm{ and }R(a,b)=\sum_{\substack{q',\\ i\in\{0,\pm 1\}} }\Delta_{q'} a\Delta_{q'+i} b.
\end{displaymath}
 Then we have

\begin{equation}
\Delta_q(ab)=\Delta_q T_a b+\Delta_q T_b a+\Delta_q R(a,b)=\Delta_q T_a b+\Delta_q \tilde R(a,b)
\end{equation} where $\tilde R(a,b)=T_b a+R(a,b)=\Sigma_{q'} S_{q'+2}b\Delta_{q'}a$. Moreover:

\begin{eqnarray}
\Delta_q (ab)=\Sigma_{|q'-q|\le 5}\Delta_q (S_{q'-1}a\Delta_{q'}b)+\Sigma_{q'> q-5}\Delta_q(S_{q'+2}b\Delta_{q'}a)\nonumber\\=\Sigma_{|q'-q|\le 5}[\Delta_q,S_{q'-1}a]\Delta_{q'}b+\Sigma_{|q'-q|\le 5}S_{q'-1}a\Delta_q\Delta_{q'}b+\Sigma_{q'> q-5}\Delta_q(S_{q'+2}b\Delta_{q'}a)\nonumber\\
=\Sigma_{|q'-q|\le 5}[\Delta_{q},S_{q'-1} a]\Delta_{q'}b+\Sigma_{|q'-q|\le 5}(S_{q'-1} a-S_{q-1}a)\Delta_q\Delta_{q'}b\nonumber\\
+\Sigma_{q'> q-5}\Delta_q (S_{q'+2} b\Delta_{q'}a)+\underbrace{\Sigma_{|q'-q|\le 5}S_{q-1}a\Delta_q\Delta_{q'}b}_{=S_{q-1}a\Delta_q b}
\label{bonydecomp}
\end{eqnarray}

\par In terms of this decomposition we can express the Sobolev norm  of an element $u$ in the space $H^s$ as:

\begin{displaymath}
\|u\|_{H^s}=\big(\|S_0 u\|_{L^2}^2+\sum_{q\in\mathbb{N}}2^{2qs}\|\Delta_q u\|_{L^2}^2\big)^{1/2}
\end{displaymath}

\par We will use the following well-known estimates:

\begin{lemma} (\cite{chemin},\cite{chemin&masmoudi})
\par {\bf (i)} (Bernstein inequalities) $$2^{-q}\|\nabla S_q u\|_{L^p}\le C \|u\|_{L^p}, \forall 1\le p\le \infty$$
$$\| \Delta_q u\|_{L^p}\le C2^{-q}\|\Delta_q \nabla u\|_{L^p}\le C\|\Delta_q u\|_{L^p}, \forall 1\le p\le \infty$$

\par {\bf (ii)}(commutator estimate) $$\| [S_{q'-1}a,\Delta_q] b \|_{L^2}\le C 2^{-q} \| \nabla  S_{q'-1}a \|_{L^\infty} \|b \|_{L^2}$$
\end{lemma}

\subsection{Proof of theorem ~\ref{theorem:reg}}

\smallskip\par {\it Step 1. Estimates of the high frequencies}
\par Applying $\Delta_q$ to the first equation in (\ref{system}) we get:

\begin{eqnarray}
\partial_t \Delta_q Q_{\alpha\beta}- \Gamma L\Delta\Delta_q Q_{\alpha\beta}-\Delta_q \Omega_{\alpha\gamma} S_{q-1} Q_{\gamma\beta}+S_{q-1}Q_{\alpha\gamma}\Delta_q\Omega_{\gamma\beta}=-\Delta_q( u_\gamma Q_{\alpha\beta,\gamma})\nonumber\\+\Gamma\Delta_q [-aQ_{\alpha\beta}+b\left(Q_{\alpha\gamma}Q_{\gamma\beta}-\frac{\delta_{\alpha\beta}}{2}\textrm{tr}(Q^2)\right)-cQ_{\alpha\beta}\textrm{tr}(Q^2)]\nonumber\\
+\Sigma_{|q'-q|\le 5} [\Delta_q;S_{q'-1}Q_{\gamma\beta}]\Delta_{q'}\Omega_{\alpha\gamma}+\Sigma_{|q'-q|\le 5}(S_{q'-1}Q_{\gamma\beta}-S_{q-1}Q_{\gamma\beta})\Delta_q\Delta_{q'}\Omega_{\alpha\gamma}+
\Sigma_{q'>q-5} \Delta_q\left(S_{q'+2}\Omega_{\alpha\gamma} \Delta_{q'}Q_{\gamma\beta}\right)\nonumber\\
-\Sigma_{|q'-q|\le 5} [\Delta_q;S_{q'-1}Q_{\alpha\gamma}]\Delta_{q'}\Omega_{\gamma\beta}-\Sigma_{|q'-q|\le 5}(S_{q'-1}Q_{\alpha\gamma}-S_{q-1}Q_{\alpha\gamma})\Delta_q\Delta_{q'}\Omega_{\gamma\beta}-
\Sigma_{q'>q-5} \Delta_q\left(S_{q'+2}\Omega_{\gamma\beta}\Delta_{q'}Q_{\alpha\gamma}\right)
\end{eqnarray}

\par Multiplying the previous equation by $-L\Delta \Delta_q Q_{\alpha\beta}$ and integrating over $\mathbb{R}^2$ and by parts we obtain:

\begin{eqnarray}
\frac{L}{2}\partial_t \|\nabla \Delta_q Q\|_{L^2}^2+\Gamma L^2\|\Delta\Delta_q Q\|_{L^2}^2+L\int \Delta_q\Omega_{\alpha\gamma} S_{q-1} Q_{\gamma\beta}\Delta\Delta_q Q_{\alpha\beta}-L\int S_{q-1} Q_{\alpha\gamma}\Delta_q \Omega_{\gamma\beta}\Delta\Delta_q Q_{\alpha\beta}\nonumber\\ =L\underbrace{\left(\Delta_q(u\nabla Q_{\alpha\beta}),\Delta\Delta_q Q_{\alpha\beta}\right)}_{\stackrel{def}{=}\mathcal{I}_1}-L\underbrace{\Sigma_{|q'-q|\le 5} \left([\Delta_q;S_{q'-1}Q_{\gamma\beta}]\Delta_{q'}\Omega_{\alpha\gamma},\Delta\Delta_q Q_{\alpha\beta}\right)}_{\stackrel{def}{=}\mathcal{I}_2}\nonumber\\-L\underbrace{\Sigma_{|q'-q|\le 5}\left((S_{q'-1}Q_{\gamma\beta}-S_{q-1}Q_{\gamma\beta})\Delta_q\Delta_{q'}\Omega_{\alpha\gamma},\Delta\Delta_q Q_{\alpha\beta}\right)}_{\mathcal{I}_3}\nonumber\\
-L\underbrace{\Sigma_{q'>q-5} \left( \Delta_q\left(S_{q'+2}\Omega_{\alpha\gamma} \Delta_{q'}Q_{\gamma\beta}\right),\Delta\Delta_q Q_{\alpha\beta}\right)}_{\stackrel{def}{=}\mathcal{I}_4}+L\underbrace{\Sigma_{|q'-q|\le 5}\left( [\Delta_q;S_{q'-1}Q_{\alpha\gamma}]\Delta_{q'}\Omega_{\gamma\beta},\Delta\Delta_q Q_{\alpha\beta}\right)}_{\stackrel{def}{=}\mathcal{I}_5}\nonumber\\
+L\underbrace{\Sigma_{|q'-q|\le 5}\left((S_{q'-1}Q_{\alpha\gamma}-S_{q-1}Q_{\alpha\gamma})\Delta_q\Delta_{q'}\Omega_{\gamma\beta},\Delta\Delta_q Q_{\alpha\beta}\right)}_{\stackrel{def}{=}\mathcal{I}_6}+L\underbrace{\Sigma_{q'>q-5}\left(  \Delta_q\left(S_{q'+2}\Omega_{\gamma\beta}\Delta_{q'}Q_{\alpha\gamma})\right),\Delta\Delta_q Q_{\alpha\beta}\right)}_{\stackrel{def}{=}\mathcal{I}_7}\nonumber\\
-L\Gamma\underbrace{\Big(\Delta_q [-aQ_{\alpha\beta}+bQ_{\alpha\gamma}Q_{\gamma\beta}-cQ_{\alpha\beta}\textrm{tr}(Q^2)],\Delta\Delta_q Q_{\alpha\beta}\Big)}_{\mathcal{I}_8}
\label{firsthalf}
\end{eqnarray}

\par Applying $\Delta_q$ to the second equation in (\ref{system}) we get:

\begin{eqnarray}
\partial_t\Delta_q u_\alpha-\nu\Delta\Delta_q u_\alpha=\partial_\alpha\Delta_q  p+L\partial_\beta\left(S_{q-1}Q_{\alpha\gamma}\Delta_q\Delta Q_{\gamma\beta}-\Delta_q\Delta Q_{\alpha\gamma} S_{q-1}Q_{\gamma\beta}\right)\nonumber\\
-L\partial_\beta\Delta_q\left(\partial_\alpha Q_{\gamma\delta}\partial_\beta Q_{\gamma\delta}-\frac{\delta_{\alpha\beta}}{3}\partial_\lambda Q_{\gamma\delta}\partial_\lambda Q_{\gamma\delta}\right)-\Delta_q(u_\beta\partial_\beta u_\alpha)\nonumber\\
+L\partial_\beta\Big(\Sigma_{|q'-q|\le 5}[\Delta_q;S_{q'-1} Q_{\alpha\gamma}]\Delta_{q'}\Delta Q_{\gamma\beta}+\Sigma_{|q'-q|\le 5}(S_{q'-1}Q_{\alpha\gamma}-S_{q-1}Q_{\alpha\gamma})\Delta_q\Delta_{q'} \Delta Q_{\gamma\beta}\Big)\nonumber\\
+L\partial_\beta\Big(\Sigma_{q'> q-5}\Delta_q(S_{q'+2}\Delta Q_{\gamma\beta}\Delta_{q'} Q_{\alpha\gamma})
-\Sigma_{|q'-q|\le 5}[\Delta_q;S_{q'-1}  Q_{\gamma\beta}]\Delta_{q'}\Delta Q_{\alpha\gamma}\Big)\nonumber\\
-L\partial_\beta\Big(\Sigma_{|q'-q|\le 5}(S_{q'-1} Q_{\gamma\beta}-S_{q-1} Q_{\gamma\beta})\Delta_q\Delta_{q'}\Delta Q_{\alpha\gamma}+\Sigma_{q'> q-5}\Delta_q(S_{q'+2}\Delta Q_{\alpha\gamma}\Delta_{q'}  Q_{\gamma\beta})\Big)
\end{eqnarray}

\par We multiply the last equation by $\Delta_q u_\alpha$, integrate over $\mathbb{R}^2$ and by parts to obtain:

\begin{eqnarray}
\frac{1}{2}\partial_t\|\Delta_q u\|_{L^2}^2+\nu\|\Delta_q\nabla u\|_{L^2}^2+L\int S_{q-1}Q_{\alpha\gamma}\Delta_q\Delta Q_{\gamma\beta}\Delta_q u_{\alpha,\beta}-L\int \Delta_q\Delta Q_{\alpha\gamma} S_{q-1}Q_{\gamma\beta}\Delta_q u_{\alpha,\beta}\nonumber\\=-\underbrace{\left(\Delta_q(u_\beta\partial_\beta u_\alpha),\Delta_q u_\alpha\right)}_{\stackrel{def}{=}\mathcal{J}_1}
+L\underbrace{\int \Delta_q\left(\partial_\alpha Q_{\gamma\delta}\partial_\beta Q_{\gamma\delta}-\frac{\delta_{\alpha\beta}}{3}\partial_\lambda Q_{\gamma\delta}\partial_\lambda Q_{\gamma\delta}\right)\Delta_q u_{\alpha,\beta}}_{\stackrel{def}{=}\mathcal{J}_2}\nonumber\\
-L\underbrace{\Sigma_{|q'-q|\le 5}\int[\Delta_q;S_{q'-1} Q_{\alpha\gamma}]\Delta_{q'}\Delta Q_{\gamma\beta}\Delta_q u_{\alpha,\beta}}_{\stackrel{def}{=}\mathcal{J}_3}
-L\underbrace{\int\Sigma_{|q'-q|\le 5}(S_{q'-1}Q_{\alpha\gamma}-S_{q-1}Q_{\alpha\gamma})\Delta_q\Delta_{q'}\Delta Q_{\gamma\beta}\Delta_q u_{\alpha,\beta}}_{\stackrel{def}{=}\mathcal{J}_4}\nonumber\\-L\underbrace{\int\Sigma_{q'> q-5}\Delta_q(S_{q'+2} \Delta Q_{\gamma\beta}\Delta_{q'} Q_{\alpha\gamma})\Delta_q u_{\alpha,\beta}}_{\stackrel{def}{=}\mathcal{J}_5}+L\underbrace{\Sigma_{|q'-q|\le 5}\int [\Delta_q;S_{q'-1}  Q_{\gamma\beta}]\Delta_{q'} \Delta Q_{\alpha\gamma}\Delta_q u_{\alpha,\beta}}_{\stackrel{def}{=}\mathcal{J}_6}\nonumber\\
+L\underbrace{\int \Sigma_{|q'-q|\le 5}(S_{q'-1} Q_{\gamma\beta}-S_{q-1} Q_{\gamma\beta})\Delta_q\Delta_{q'}\Delta Q_{\alpha\gamma}\Delta_q u_{\alpha,\beta}}_{\stackrel{def}{=}\mathcal{J}_7}+L\underbrace{\int \Sigma_{q'>q-5}\Delta_q(S_{q'+2}\Delta Q_{\alpha\gamma}\Delta_{q'} Q_{\gamma\beta})\Delta_q u_{\alpha,\beta}}_{\stackrel{def}{=}\mathcal{J}_8}
\label{secondhalf}
\end{eqnarray}

\par Summing (\ref{firsthalf}) and (\ref{secondhalf}) and using Lemma ~\ref{anulare} we get:

\begin{eqnarray}
\partial_t \left(\frac{L}{2}\|\nabla \Delta_q Q\|_{L^2}^2+\frac{1}{2}\|\Delta_q u\|_{L^2}^2\right)+\nu\|\Delta_q\nabla u\|_{L^2}^2+\Gamma L^2\|\Delta\Delta_q Q\|_{L^2}^2=\sum_{i=1}^{8} \mathcal{I}_i+\sum_{j=1}^8\mathcal{J}_j
\end{eqnarray} 

\par We denote by $\varphi(t)\stackrel{def}{=}L\|\nabla Q\|_{H^s}^2+\|u\|_{H^s}^2$ with $\varphi_1(t)\stackrel{def}{=}L\|S_0 \nabla Q\|_{L^2}^2+\|S_0 u\|_{L^2}^2$ the low-frequency part of $\varphi$ and  $\varphi_2(t)\stackrel{def}{=}\varphi(t)-\varphi_1(t)$ the high-frequency part of $\varphi$.

\par The last equality leads to the following estimate, whose technical proof is postponed to the Appendix $A$:

\begin{eqnarray}
\frac{1}{2}\frac{d}{dt}\varphi_2+\sum_{q\in\mathbb{N}} 2^{2qs}\Big(\frac{\Gamma L^2}{2}\|\Delta\Delta_q Q\|_{L^2}^2+\frac{\nu}{2}\|\nabla\Delta_q u\|_{L^2}^2\Big)\nonumber\\ 
\leq C(1+\|\nabla Q\|_{L^\infty}^2+\|u\|_{L^\infty}^2+\|Q\|_{L^\infty}^2\Big))(\|\nabla Q\|_{H^s}^2+\|u\|_{H^s}^2)+\frac{\Gamma L^2}{50}\|\Delta Q\|_{H^s}^2+\frac{\nu}{50}\|\nabla u\|_{H^s}^2
\label{longest}
\end{eqnarray}

\bigskip\par{\it Step 2. Estimates of the low frequencies}
\par This is much easier than the previous step. We apply $S_0$ to the first equation in (\ref{system}), multiply by $-LS_0\Delta Q_{\alpha\beta}$, take the trace, integrate over $\mathbb{R}^2$ and by parts  and we get:
\begin{eqnarray}
\frac{L}{2}\partial_t \|S_0\nabla  Q\|_{L^2}^2+\Gamma L^2\|\Delta S_0 Q\|_{L^2}^2\le C\|u\|_{L^\infty}\|S_0\nabla Q\|_{L^2}\|\Delta S_0 Q\|_{L^2}+C\|Q\|_{L^\infty}\|\nabla u\|_{L^2}\|\Delta S_0 Q\|_{L^2}\nonumber\\
+C\|\nabla S_0 Q\|_{L^2}^2\big(1+\|Q\|_{L^\infty}+\|Q\|_{L^\infty}^2\big)\nonumber
\end{eqnarray} hence

\begin{equation}
\frac{L}{2}\partial_t \|S_0\nabla Q\|_{L^2}^2+\frac{\Gamma L^2}{2}\|\Delta S_0 Q\|_{L^2}^2\le C\|S_0\nabla Q\|_{L^2}^2\Big(\|u\|_{L^\infty}^2+1+\|Q\|_{L^\infty}+\|Q\|_{L^\infty}^2\Big)+C\|u\|_{H^s}^2\|Q\|_{L^\infty}^2
\label{eq:low1}
\end{equation}

\par We aply $S_0$ to the second equation in (\ref{system}), multiply by $S_0 u$ and integrate over $\mathbb{R}^2$ and by parts to obtain:
\begin{displaymath}
\frac{1}{2}\partial_t \|S_0 u\|_{L^2}^2+\nu\|\nabla S_0 u\|_{L^2}^2\le C\|u\|_{L^\infty}\|\nabla u\|_{L^2}\|S_0 u\|_{L^2}+C\|\nabla S_0 u\|_{L^2}
\Big(\|\nabla Q\|_{L^\infty}\|\nabla Q\|_{L^2}+\|Q\|_{L^\infty}\|\Delta Q\|_{L^2}\Big)
\end{displaymath} hence

\begin{equation}
\frac{1}{2}\partial_t\|S_0 u\|_{L^2}^2+\frac{\nu}{2}\|\nabla S_0 u\|_{L^2}^2\le C\|u\|_{H^s}^2\Big(1+\|u\|_{L^\infty}\Big)+C\|\nabla Q\|_{H^s}^2\Big(\|\nabla Q\|_{L^2}^2+\|Q\|_{L^\infty}^2\Big)
\label{eq:low2}
\end{equation}

\par Summing (\ref{eq:low1}) and (\ref{eq:low2}) we obtain:
\begin{equation}
\partial_t \varphi_1+\frac{\nu}{2}\|\nabla S_0 u\|_{L^2}^2+\frac{\Gamma L^2}{2}\|\Delta S_0 Q\|_{L^2}^2\le \|u\|_{L^\infty}^2\varphi+C\Big(1+\|Q\|_{L^\infty}^2+\|\nabla Q\|_{L^2}^2\Big)\varphi
\label{est:low}
\end{equation}

\medskip\par{\it Step 3. The estimates of the high norms}

\par Summing (\ref{longest}) and (\ref{est:low}) we obtain:

$$\frac{1}{2}\varphi'(t)\nonumber\\ 
\leq C(\|\nabla Q\|_{L^\infty}^2+\|u\|_{L^\infty}^2)\varphi+C\Big(1+\|Q\|_{L^\infty}^2+\|\nabla Q\|_{L^2}^2\Big)\varphi
$$

\par Now we use a fundamental ingredient in the global existence, namely the logarithmic estimate (see \cite{brezisgallouet}), for $s>1$,
$$\|\nabla Q\|_{L^\infty}+\|u\|_{L^\infty}\leq C(\|Q\|_{H^2}+\|u\|_{H^1})\sqrt{\ln(e+\frac{\|\nabla Q\|_{H^s}^2+\|u\|_{H^s}^2}{\|Q\|_{H^2}+\|u\|_{H^1}})},$$
and be denoting $f(t)\stackrel{def}{=}\|Q\|_{H^2}^2+\| u\|^2_{H^1}$ and$g(t)\stackrel{def}{=}1+\|Q\|_{L^\infty}^2+\|\nabla Q\|_{L^2}^2$ we obtain
$$\varphi'(t)\leq C f(t)\big(\ln(e+\frac{\varphi(t)}{f(t)})+g(t)\big)\varphi(t).$$

\par Observing that the function $h(x)\stackrel{def}{=}x\ln(e+\frac{\varphi}{x})$ is increasing the last relation implies:
$$\varphi'(t)\le C(1+f(t))\Big(\ln(e+\varphi(t))+g(t)\Big)\varphi(t)$$

By integrating this differential inequality, we obtain:
$$\varphi(t)\leq (e+\|Q_0\|_{H^s}^2+\|u_0\|_{H^s}^2)^{\Big(e^{\int_0^t (1+f(\tau))d\tau}C\big(1+\int_0^t g(\tau)\,d\tau\big)\Big)}\leq C(t,\bar Q, \bar u,s),$$
and this uniform bound, imply the global existence of a regular solution for regular enough initial data. Taking into account Proposition ~\ref{prop:aprioriest} we have that $\int_0^t f(s)ds$ increases exponentially and this gives the rate in (\ref{rate:3exp}). $\Box$

\section{Weak-Strong uniqueness in 2D}
In this section we consider a global weak solution and a strong one, starting from the same  initial data $(\bar Q,\bar u)\in H^{s+1}(\mathbb{R}^2)\times H^s(\mathbb{R}^2)$ with $s>1$ and we show that they are the same. More precisely:

\begin{proposition} Let $(\bar Q,\bar u)\in H^{s+1}(\mathbb{R}^2)\times H^s(\mathbb{R}^2)$ with $s>1$. By Proposition ~\ref{prop:weak} there exists a weak solution $(Q_1, u_1)$ of the system (\ref{system}), subject to restriction (\ref{c+}) and starting from initial data $(\bar Q,\bar u)$, such that 

 \begin{equation}Q_1\in L^\infty_{loc}(\mathbb{R}_+;H^1(\mathbb{R}^2))\cap L^2_{loc}(\mathbb{R}_+;  H^2(\mathbb{R}^2))\textrm{ and } u_1\in L^\infty_{loc}(\mathbb{R}_+;L^2(\mathbb{R}^2))\cap L^2_{loc}(\mathbb{R}_+; H^1(\mathbb{R}^2))
 \label{sol1}
 \end{equation}

\par  Theorem ~\ref{theorem:reg} gives the existence of a strong solution $(Q_2,u_2)$ such that

 \begin{equation}Q_2\in L^\infty_{loc}(\mathbb{R}_+; H^{s+1}(\mathbb{R}^2)\cap L^2_{loc}(\mathbb{R}_+; H^{s+2}(\mathbb{R}^2))\textrm{ and }u_2\in L^\infty(\mathbb{R}_+; H^s(\mathbb{R}^2))\cap L^2(\mathbb{R}_+; H^{s+1}(\mathbb{R}^2))
 \label{sol2}
 \end{equation} with $s>1$ and the same initial data $(\bar Q,\bar u)\in H^{s+1}(\mathbb{R}^2)\times H^s(\mathbb{R}^2)$. Then  $(Q_1,u_1)=(Q_2, u_2)$. 
\end{proposition}

\smallskip\par{\bf Proof.}
We denote by $\delta Q=Q_1-Q_2$ and $\delta u=u_1-u_2$ which verify the following system
\begin{equation}
\left\{\begin{array}{l}
      (\partial_t+  \delta u\nabla )\delta Q-\delta \Omega \delta Q+\delta Q\delta \Omega+\delta u\nabla Q_2+u_2\nabla\delta Q+Q_2\delta\Omega+\delta Q\Omega_2-\delta\Omega Q_2-\Omega_2\delta Q\\=\Gamma\Big(L\Delta \delta Q-a\delta Q+b[\delta QQ_1+Q_2\delta Q-\frac{\textrm{tr}\big(\delta QQ_1+Q_2\delta Q\big)}{2}Id]-c\delta Q\textrm{tr}(Q_1^2)-cQ_2\big[\textrm{tr}(Q_1\delta Q+\delta QQ_2)\big]\Big) \\
      \partial_t \delta u +\mathcal{P}(\delta u\nabla \delta u)=\nu\Delta \delta u-L\mathcal{P}\big(\nabla\cdot(\nabla \delta Q\nabla \delta Q-\frac{1}{2}|\nabla \delta Q|^2)\big)+L\mathcal{P}\big(\nabla\cdot(\delta Q\Delta \delta Q-\Delta \delta Q\delta Q)\big)\\
 -\mathcal{P}(u_2\nabla\delta u+\delta u\nabla u_2)-L \mathcal{P}\bigg(\nabla\cdot\Big((\nabla\delta Q\nabla Q_2+\nabla Q_2\nabla\delta Q)-\frac{1}{2}\textrm{tr}\big(\nabla\delta Q\nabla Q_2+Q_2\nabla\delta Q\big)Id\Big)\bigg)\\+L\mathcal{P}\big(\nabla\cdot(\delta Q\Delta Q_2+Q_2\Delta \delta Q-\Delta\delta Q Q_2-\Delta Q_2\delta Q)\big)\\
     \end{array}\right.
      \label{approxsystemdif}
      \end{equation}
      
\smallskip  
\par We proceed similarly as in the proof of Proposition ~\ref{prop:Lyapunov}, namely we multiply  the first equation in (\ref{approxsystemdif}) to the right by $-L\Delta\delta Q+\delta Q$, integrate over $\mathbb{R}^2$ and by parts, take the trace and sum with the second equation in (\ref{approxsystemdif}) multiplied by $\delta u$ and integrated over $\mathbb{R}^2$ and by parts.
 Taking into account the cancellations analogous to the ones in (\ref{Lyapunovcancellation}) we obtain:
 \begin{eqnarray}
 \frac{d}{dt}\int_{\mathbb{R}^2}\frac{L}{2}|\nabla\delta Q(x)|^2+\frac{1}{2}|\delta Q(x)|^2+\frac{1}{2}|\delta u(x)|^2\,dx+\int_{\mathbb{R}^2}\nu|\nabla\delta u(x)|^2+\Gamma L^2 |\Delta\delta Q(x)|^2\,dx\nonumber\\
 =L\int_{\mathbb{R}^2}\textrm{tr}\Big(\big[\delta u\nabla Q_2+u_2\nabla\delta Q+\delta Q\Omega_2-\Omega_2\delta Q\big]\Delta\delta Q\Big)\,dx+\underbrace{L\int_{\mathbb{R}^2}\textrm{tr}\Big(\big[Q_2\delta\Omega-\delta\Omega Q_2\big]\Delta\delta Q\Big)\,dx}_{\mathcal{A}}\nonumber\\-a\Gamma L\int_{\mathbb{R}^2}|\nabla\delta Q(x)|^2\,dx 
-b\Gamma L\int_{\mathbb{R}^2}\textrm{tr}\Big(\big(\delta Q(x) Q_1(x) +Q_2(x)\delta Q(x)\big)\Delta\delta Q(x)\Big)\,dx\nonumber\\
+c\Gamma L\int_{\mathbb{R}^2}\textrm{tr}\big(\delta Q\Delta\delta Q\big)\textrm{tr}(Q_1^2)\,dx+c\Gamma L\int_{\mathbb{R}^2}\textrm{tr}(Q_2\Delta\delta Q)\textrm{tr}(Q_1\delta Q+\delta Q Q_2)\,dx\nonumber
\end{eqnarray}
\begin{eqnarray}
-\int_{\mathbb{R}^2}\textrm{tr}\big(\delta u\nabla Q_2\delta Q\big)\,dx-\int_{\mathbb{R}^2}\textrm{tr}\big(Q_2\delta\Omega\delta Q\big)\,dx-\underbrace{\int_{\mathbb{R}^2}\textrm{tr}\big(\delta Q\Omega_2\delta Q\big)dx}_{\mathcal{I}}\nonumber\\
+\int_{\mathbb{R}^2}\textrm{tr}\big(\delta\Omega Q_2\delta Q\big)\,dx+\underbrace{\int_{\mathbb{R}^2}\textrm{tr}\big(\Omega_2(\delta Q)^2\big)\,dx}_{\mathcal{II}}-\Gamma L\int_{\mathbb{R}^2}|\nabla Q|^2\,dx\nonumber\\
\-a\Gamma\int_{\mathbb{R}^2}|\delta Q|^2\,dx+b\Gamma \int_{\mathbb{R}^2}\textrm{tr}\big(\delta Q Q_1\delta Q+Q_2(\delta Q)^2\big)\,dx\nonumber\\
\end{eqnarray}
\begin{eqnarray}
-c\Gamma\int_{\mathbb{R}^2}\textrm{tr}(Q_1)^2|\delta Q|^2\,dx-c\Gamma\int_{\mathbb{R}^2}\textrm{tr}(Q_2\delta Q)\textrm{tr}(Q_1\delta Q+\delta QQ_2)\,dx\nonumber\\
-\int_{\mathbb{R}^2}\big(u_2\nabla \delta u+\delta u\nabla u_2)\delta u\,dx+L\int_{\mathbb{R}^2}\big(\nabla\delta Q\nabla Q_2+\nabla Q_2\nabla\delta Q\big)\cdot\nabla \delta u\,dx\nonumber\\-L\int_{\mathbb{R}^2} \Big(\delta Q\Delta Q_2-\Delta Q_2\delta Q\Big)\cdot\nabla\delta u\,dx-\underbrace{L\int_{\mathbb{R}^2} \Big(Q_2\Delta\delta Q-\Delta\delta Q Q_2\Big)\cdot\nabla\delta u\,dx}_{\mathcal{AA}}
 \end{eqnarray}

\par Let us observe that Lemma ~\ref{anulare} implies $\mathcal{A}-\mathcal{AA}=0$. Also $\mathcal{I}+\mathcal{II}=0$ and then  we easily obtain
\begin{eqnarray}
\frac 12\frac{d}{dt}(L\|\nabla\delta Q\|^2_{L^2}+\|\delta Q\|_{L^2}^2+\|\delta u\|_{L^2}^2)+\Gamma L^2\|\Delta \delta Q\|_{L^2}^2+\nu\|\nabla \delta u\|_{L^2}^2\leq L\|\Delta\delta Q\|_{L^2}\|\delta u\|_{L^2}\|\nabla Q_2\|_{L^\infty}\nonumber\\
+L\|u_2\|_{L^\infty}\|\nabla\delta Q\|_{L^2}\|\Delta\delta Q\|_{L^2}+2L\|\delta Q\|_{L^2}\|\Omega_2\|_{L^\infty}\|\Delta\delta Q\|_{L^2}\nonumber\\
|a|\Gamma L\|\nabla\delta Q\|_{L^2}^2+|b|\Gamma L\|\Delta\delta Q\|_{L^2}\|\delta Q\|_{L^4} \|Q_1\|_{L^4} +|b|\Gamma L\|Q_2\|_{L^\infty}\|\delta Q \|_{L^2}\|\Delta\delta Q\|_{L^2}\nonumber\\
+c\Gamma L\|\delta Q\|_{L^4}\|\Delta\delta Q\|_{L^2}\|Q_1\|_{L^8}^2+c\Gamma L\|Q_2\|_{L^\infty}\|\Delta\delta Q\|_{L^2}\big(\|Q_1\|_{L^4}+\|Q_2\|_{L^4}\big)\|\delta Q\|_{L^4}\nonumber\\
+\|\nabla Q_2\|_{L^\infty}\|\delta u\|_{L^2}\|\delta Q\|_{L^2}+2\|Q_2\|_{L^\infty}\|\nabla\delta u\|_{L^2}\|\delta Q\|_{L^2}\nonumber\\
+|a|\Gamma\|\delta Q\|_{L^2}^2+\Gamma\big(|b|+c\|Q_2\|_{L^\infty}\big)\big(\|Q_1\|_{L^2}+\|Q_2\|_{L^2}\big)\|\delta Q\|_{L^4}^2\nonumber\\
+\|\delta u\|_{L^2}^2\|\nabla u_2\|_{L^\infty}+2L\|\nabla Q_2\|_{L^\infty}\|\nabla\delta Q\|_{L^2}\|\nabla \delta u\|_{L^2}+2L\|\Delta Q_2\|_{L^\infty}\|\delta Q\|_{L^2}\|\nabla\delta u\|_{L^2}\nonumber
\end{eqnarray}
\begin{eqnarray}
\le\frac{\nu}{2}\|\nabla \delta u\|_{L^2}+\frac{\Gamma L^2}{2}\|\Delta\delta Q\|_{L^2}^2+C\underbrace{\Big(\|\nabla u_2\|_{L^\infty}+\|\nabla Q_2\|_{L^\infty}^2\Big)}_{\mathcal{J}_1}\|\delta u\|_{L^2}^2\nonumber\\+C\underbrace{\Big(1+\|\nabla u_2\|_{L^\infty}^2+\|\nabla Q_2\|_{L^\infty}^2+\|Q_2\|_{L^\infty}^2+\|\Delta Q_2\|_{L^\infty}^2\Big)}_{\mathcal{J}_2}\|\delta Q\|_{L^2}^2+C\underbrace{\Big(1+\|u_2\|_{L^\infty}^2+\|\nabla Q_2\|_{L^\infty}^2\Big)}_{\mathcal{J}_3}\|\nabla \delta Q\|_{L^2}^2\nonumber\\+C\underbrace{\bigg(\|Q_1\|_{L^4}^2+\|Q_1\|_{L^8}^4+\|Q_2\|_{L^\infty}^2\big(\|Q_1\|_{L^4}^2+\|Q_2\|_{L^4}^2\big)+\Gamma\big(|b|+c\|Q_2\|_{L^\infty}\big)\big(\|Q_1\|_{L^2}+\|Q_2\|_{L^2}\big)\bigg)}_{\mathcal{J}_4}\|\delta Q\|_{L^4}^2
\end{eqnarray}

\par We are in $2D$ so $\|\delta Q\|_{L^4}^2$ is controlled by $\|\delta Q\|_{L^2}^2+\|\nabla \delta Q\|_{L^2}^2$. The hypothesis, namely relations (\ref{sol1}) and (\ref{sol2}), ensure that the terms $\mathcal{J}_i,i=1,2,3,4$ are integrable in time  thus using the last inequality and  Gronwall Lemma we obtain the uniqueness of the solution. $\Box$

\bigskip\noindent{\bf Acknowledgements} MP and AZ thank John M. Ball for stimulating discussions. MP  gratefully acknowledges the hospitality of Oxford University's OXPDE Center. AZ acknowledges the support of  the EPSRC Science and Innovation award to the Oxford Center for Nonlinear PDE (EP/E035027/1).

\appendix
\section{Proof of estimate (\ref{longest})}
\par In the following, $a_q(t)$  denotes a sequence in $l_{q}^2$  for all $t>0$ and $b_q(t)$ is a sequence in  $l^1_q$, $\forall t\ge 0$, sequences that can change from one line to the next. Moreover $\|\big(a_q(t)\big)_{q\in\mathbb{N}}\|_{l^2},\|\big(b_q(t)\big)_{q\in\mathbb{N}}\|_{l^1}\le C$ where the constant $C$ is independent of $t\ge 0$.
\begin{eqnarray}
|\mathcal{I}_1|=|\left(\Delta_q(u\nabla Q_{\alpha\beta}),\Delta_q\Delta Q_{\alpha\beta}\right)|\stackrel{(\ref{bonydecomp})}{=}|\underbrace{\int S_{q-1} u \Delta_q\nabla Q_{\alpha\beta}\Delta_q\Delta Q_{\alpha\beta}}_{\stackrel{def}{=}\mathcal{I}_{1a}}
+\underbrace{\sum_{|q'-q|\le 5}\left([\Delta_q;S_{q'-1} u]\Delta_{q'}\nabla Q_{\alpha\beta},\Delta_q\Delta Q_{\alpha\beta}\right)}_{\stackrel{def}{=}\mathcal{I}_{1b}}\nonumber\\
+\underbrace{\sum_{|q'-q|\le 5}\left((S_{q'-1} u-S_{q-1} u)\Delta_q\Delta_{q'}\nabla Q_{\alpha\beta},\Delta\Delta_q Q_{\alpha\beta}\right)}_{\stackrel{def}{=}\mathcal{I}_{1c}}+\underbrace{\sum_{q'\ge q-5}(\Delta_q(S_{q'+2}\nabla Q_{\alpha\beta}\Delta_{q'} u),\Delta_q\Delta Q_{\alpha\beta})}_{\mathcal{I}_{1d}}|
\end{eqnarray}

\begin{eqnarray}
|\mathcal{I}_{1a}|\le C \|u\|_{L^\infty}\|\Delta_q\nabla Q\|_{L^2}\|\Delta\Delta_q Q\|_{L^2}\le C2^{-2qs}b_q(t) \|u\|_{L^\infty}\|\nabla Q\|_{H^s}\|\Delta Q\|_{H^s}
\end{eqnarray}

\begin{eqnarray}
|\mathcal{I}_{1b}|\le \sum_{|q'-q|\le 5}\| [\Delta_q;S_{q'-1} u]\Delta_{q'}\nabla Q_{\alpha\beta}\|_{L^2}\|\Delta_q \Delta Q_{\alpha\beta}\|_{L^2}\le \sum_{|q'-q|\le 5} 2^{-q}\|\nabla S_{q'-1}u\|_{L^\infty}\|\nabla\Delta_{q'} Q_{\alpha\beta}\|_{L^2}\|\Delta_q\Delta Q_{\alpha\beta}\|_{L^2}\nonumber\\
\le C\|u\|_{L^\infty} 2^{-2qs}b_q(t)\|\nabla Q\|_{H^s}\|\Delta Q\|_{H^s}
\end{eqnarray}

\begin{eqnarray}
|\mathcal{I}_{1c}|\le C\|u\|_{L^\infty}\|\Delta_q\nabla Q\|_{L^2}\|\Delta\Delta_q Q\|_{L^2}\le C2^{-2qs}b_q(t) \|u\|_{L^\infty}\|\nabla Q\|_{H^s}\|\Delta Q\|_{H^s}
\end{eqnarray}

\begin{eqnarray}
|\mathcal{I}_{1d}|\le\sum_{q'> q-5}|(\Delta_q(S_{q'+2}\nabla Q_{\alpha\beta}\Delta_{q'} u),\Delta_q\Delta Q_{\alpha\beta})|
\le \|\nabla Q\|_{L^\infty}\sum_{q'> q-5}2^{-(q'+q)s}2^{q's}\|\Delta_{q'}u\|_{L^2}2^{qs}\|\Delta_q\Delta Q\|_{L^2}\nonumber\\
\le\|\nabla Q\|_{L^\infty} \sum_{q'> q-5}2^{-(q'+q)s}b_{q'}(t)\|u\|_{H^s}\|\Delta Q\|_{H^s}\le
 C\|\nabla Q\|_{L^\infty}2^{-2qs}\tilde b_q(t)\|u\|_{H^s}\|\Delta Q\|_{H^s}
\end{eqnarray} where $\tilde b_q(t)=\sum_{q'\ge q-5}2^{-(q'-q)s}b_{q'}(t)$.

\begin{eqnarray}
|\mathcal{I}_2|=|\sum_{|q'-q|\le 5} \left([\Delta_q;S_{q'-1}Q_{\gamma\beta}]\Delta_{q'}\Omega_{\alpha\gamma},\Delta\Delta_q Q_{\alpha\beta}\right)|\le \sum_{|q'-q|\le 5}  2^{-q}\|S_{q'-1}\nabla Q_{\gamma\beta}\|_{L^\infty}\|\Delta_{q'} \Omega_{\alpha\gamma}\|_{L^2}\|\Delta\Delta_q Q_{\alpha\beta}\|_{L^2}\nonumber\\
\le \sum_{|q'-q|\le 5}C 2^{-q}\|S_{q'-1}\nabla Q_{\gamma\beta}\|_{L^\infty}2^{q'}\|\Delta_{q'}u\|_{L^2}\|\Delta\Delta_q Q_{\alpha\beta}\|_{L^2}\le  C\sum_{|q'-q|\le 5}\|\nabla Q\|_{L^\infty}\|\Delta_{q'} u\|_{L^2}\|\Delta\Delta_q Q\|_{L^2}\nonumber\\
\le  C 2^{-2qs} b_q(t)\|\nabla Q\|_{L^\infty}\|u\|_{H^s}\|\Delta Q\|_{H^s}
\end{eqnarray}

\begin{eqnarray}
|\mathcal{I}_3|=|\sum_{|q'-q|\le 5}\left((S_{q'-1}Q_{\gamma\beta}-S_{q-1}Q_{\gamma\beta})\Delta_q\Delta_{q'}\Omega_{\alpha\gamma}, \Delta\Delta_q Q_{\alpha\beta}\right)|\nonumber\\
\le \sum_{|q'-q|\le 5} \|\left(S_{q'-1}Q_{\gamma\beta}-S_{q-1}Q_{\gamma\beta}\right)\Delta_q\Delta_{q'}\Omega_{\alpha\gamma}\|_{L^2}\| \Delta\Delta_q Q_{\alpha\beta}\|_{L^2}\nonumber\\
\le C\sum_{|q'-q|\le 5}\|S_{q'-1}Q_{\gamma\beta}-S_{q-1}Q_{\gamma\beta}\|_{L^\infty}\|\Delta_q\Omega_{\alpha\gamma}\|_{L^2}
\| \Delta\Delta_q Q_{\alpha\beta}\|_{L^2}\nonumber\\
\le \sum_{|q'-q|\le 5}2^{-q'}\|\tilde\Delta_{q'}\nabla Q_{\gamma\beta}\|_{L^\infty}\|\Delta_q\Omega_{\alpha\gamma}\|_{L^2}\| \Delta\Delta_q Q_{\alpha\beta}\|_{L^2}\nonumber\\
\le C2^{-q}\|\nabla Q\|_{L^\infty}2^q\|\Delta_q u\|_{L^2}\|\Delta_q\Delta Q\|_{L^2}\le C2^{-2qs}b_q(t)\|\nabla Q\|_{L^\infty}\|u\|_{H^s}\|\Delta Q\|_{H^s}
\end{eqnarray} where $\tilde\Delta_q=\sum_{|i|\le 5}\Delta_q$.

\begin{eqnarray}
|\mathcal{I}_4|=|\sum_{q'> q-5}\left(\Delta_q\left(S_{q'+2}\Omega_{\alpha\gamma}\Delta_{q'}Q_{\gamma\beta}\right),\Delta_q \Delta Q_{\alpha\beta}\right)|\le \sum_{q'> q-5}\|\Delta_q(S_{q'+2}\Omega_{\alpha\gamma}\Delta_{q'}Q_{\gamma\beta})\|_{L^2}\|\Delta_q\Delta Q_{\alpha\beta}\|_{L^2}\nonumber\\
\le \sum_{q'> q-5}\|S_{q'+2}\Omega_{\alpha\gamma}\|_{L^\infty}\|\Delta_{q'}Q_{\gamma\beta}\|_{L^2}\|\Delta_q\Delta Q_{\alpha\beta}\|_{L^2}\le  \sum_{q'> q-5} C2^{q'}\|S_{q'+2} u\|_{L^\infty}\|\Delta_{q'}Q_{\gamma\beta}\|_{L^2} \|\Delta_q \Delta Q_{\alpha\beta}\|_{L^2}\nonumber\\
\le \|u\|_{L^\infty}C\sum_{q'> q-5} 2^{-2q's}2^{q'(s+1)}\|\Delta_{q'}Q_{\gamma\beta}\|_{L^2}2^{q's}\|\Delta_q\Delta Q_{\alpha\beta}\|_{L^2}\le C\|u\|_{L^\infty}\|\nabla Q\|_{H^s}\|\Delta Q\|_{H^s}\Sigma_{q'\ge q-5} 2^{-2q's}a_{q'}(t)\bar a_{q'}(t)\nonumber\\
\le C2^{-2qs} b_q(t)\|u\|_{L^\infty}\|\nabla Q\|_{H^s}\|\Delta Q\|_{H^s}\nonumber
\end{eqnarray} where  $b_q(t)=\sum_{q'>q-5} 2^{-2(q'-q)s}a_{q'}(t)\bar a_{q'}(t)$.

\smallskip\par The term $\mathcal{I}_k, k=5,6,7$ is estimated exactly as the term $\mathcal{I}_{k-3}$ that we have  already studied above.

\begin{eqnarray}
|\mathcal{I}_8|=|\Big(\Delta_q\big(a\nabla Q_{\alpha\beta}-b\nabla Q_{\alpha\gamma}Q_{\gamma\beta}-bQ_{\alpha\gamma}\nabla Q_{\gamma\beta}+c\nabla Q_{\alpha\beta}\textrm{tr}(Q^2)+2cQ_{\alpha\beta}\nabla Q_{\gamma\delta}Q_{\gamma\delta}\big),\nabla \Delta_q Q_{\alpha\beta}\Big)|\nonumber\\
\le \Big(|a|+2|b|\|Q\|_{L^\infty}+3c\|Q\|_{L^\infty}^2\Big)\|\nabla\Delta_q Q\|_{L^2}^2\le \Big(|a|+2|b|\|Q\|_{L^\infty}+3c\|Q\|_{L^\infty}^2\Big) 2^{-2qs}b_q(t)\|\nabla Q\|_{H^s}^2
\end{eqnarray}

\begin{eqnarray}
|\mathcal{J}_1|=|(\Delta_q(u\nabla u),\Delta_q u)|=\underbrace{|\int S_{q-1}u\nabla\Delta_q u\cdot\Delta_q u|}_{\mathcal{J}_{1a}}+\underbrace{\sum_{|q'-q|\le 5} |\int [\Delta_q;S_{q'-1} u]\Delta_{q'}\nabla u \Delta_q u|}_{\mathcal{J}_{1b}}\nonumber\\
+\underbrace{\sum_{|q'-q|\le 5} |\int(S_{q'-1}u-S_{q-1}u)\Delta_q\Delta_{q'}\nabla u\Delta_q u|}_{\mathcal{J}_{1c}}+\underbrace{\sum_{q'> q-5}|\int \Delta_q (S_{q'+2}\nabla u\cdot\Delta_{q'}u)\Delta_q u|}_{\mathcal{J}_{1d}}
 \end{eqnarray}

with
\begin{eqnarray}
|\mathcal{J}_{1a}|\le \|S_{q-1}u\|_{L^\infty}\|\Delta_q\nabla u\|_{L^2}\|\Delta_q u\|_{L^2}\le \|u\|_{L^\infty} 2^{-2qs} b_q(t)\|\nabla u\|_{H^s}\|u\|_{H^s}
\end{eqnarray}

\begin{eqnarray}
|\mathcal{J}_{1b}|=|\sum_{|q'-q|\le 5} \int [\Delta_q;S_{q'-1} u]\Delta_{q'}\nabla u \Delta_q u|\le C 2^{-q}\|S_{q-1}\nabla u\|_{L^\infty}\|\Delta_{q'}\nabla u\|_{L^2}\|\Delta_q u\|_{L^2}\nonumber\\
\le C \|u\|_{L^\infty}2^{-2qs}b_q(t)\|\nabla u\|_{H^s}\|u\|_{H^s}
\end{eqnarray}

\begin{eqnarray}
|\mathcal{J}_{1c}|\le \sum_{|q'-q|\le 5}\|(S_{q'-1}-S_{q-1})u\|_{L^\infty}\|\Delta_q\nabla u\|_{L^2}\|\Delta_q u\|_{L^2}\le C \|u\|_{L^\infty} 2^{-2qs} b_q(t)\|\nabla u\|_{H^s}\|u\|_{H^s}
\end{eqnarray}

\begin{eqnarray}
|\mathcal{J}_{1d}|=|\sum_{q'> q-5}\left(\Delta_q\left(S_{q'+2}\nabla u\Delta_{q'}u\right),\Delta_q u\right)|\le \sum_{q'> q-5}\|\Delta_q(S_{q'+2}\nabla u\Delta_{q'}u)\|_{L^2}\|\Delta_q u\|_{L^2}\nonumber\\
\le \sum_{q'> q-5}\|S_{q'+2}\nabla u\|_{L^\infty}\|\Delta_{q'}u\|_{L^2}\|\Delta_q u\|_{L^2}\le  \sum_{q'> q-5} C2^{q'}\|S_{q'+2} u\|_{L^\infty}\|\Delta_{q'}u\|_{L^2} \|\Delta_q u\|_{L^2}\nonumber\\
\le C\|u\|_{L^\infty}\sum_{q'> q-5} c2^{-2q's}2^{q's}\|\Delta_{q'}\nabla u\|_{L^2}2^{q's}\|\Delta_q u\|_{L^2}\nonumber\\
\le C\|u\|_{L^\infty}\|\nabla u\|_{H^s}\|u\|_{H^s}\sum_{q'> q-5} c2^{-2q's}a_{q'}(t)\bar a_{q'}(t)
\le C2^{-2qs} b_q(t)\|u\|_{L^\infty}\|\nabla u\|_{H^s}\|u\|_{H^s}
\end{eqnarray} where $b_q(t)=\sum_{q'> q-5} 2^{-2(q'-q)s}a_{q'}(t)\bar a_{q'}(t)\in l^1_q, \forall t\ge 0$.

\begin{eqnarray}
|\mathcal{J}_2|=|\int \Delta_q\left(\partial_\alpha Q_{\gamma\delta}\partial_\beta Q_{\gamma\delta}\right)\Delta_q u_{\alpha,\beta}|\le \|\Delta_q\left(\partial_\alpha Q_{\gamma\delta}\partial_\beta Q_{\gamma\delta}\right)\|_{L^2}\|\Delta_q\nabla u\|_{L^2}\nonumber\\
\le C2^{-2qs}b_q(t)\|\partial_\alpha Q_{\gamma\delta}\partial_\beta Q_{\gamma\delta}\|_{H^s}\|\nabla u\|_{H^s}
\le C2^{-2qs}b_q(t)\|\nabla Q\|_{L^\infty}\|\nabla Q\|_{H^s}\|\nabla u\|_{H^s}
\end{eqnarray}

\begin{eqnarray}
|\mathcal{J}_3|=|\sum_{|q'-q|\le 5}\int[\Delta_q;S_{q'-1} Q_{\alpha\gamma}]\Delta_{q'}\Delta Q_{\gamma\beta}\Delta_q u_{\alpha,\beta}|\le \sum_{|q'-q|\le 5} \|[\Delta_q;S_{q'-1}Q_{\alpha\gamma}]\Delta_{q'}\Delta Q_{\gamma\beta}\|_{L^2}\|\Delta_q\nabla u\|_{L^2}\nonumber\\
\le C2^{-q}\|S_{q'-1}\nabla Q_{\alpha\gamma}\|_{L^\infty}\|\Delta_q\Delta Q_{\gamma\beta}\|_{L^2}\|\Delta_q\nabla u\|_{L^2}\le C\|\nabla Q\|_{L^\infty}\|\Delta_q\nabla Q\|_{L^2}\|\Delta_q\nabla u\|_{L^2}\nonumber\\
\le C2^{-2qs}b_q(t)\|\nabla Q\|_{L^\infty}\|\nabla Q\|_{H^s}\|\nabla u\|_{H^s}
\end{eqnarray}

Concerning the term $\mathcal{J}_4$ we use that  $(S_{q'-1}Q_{\alpha\gamma}-S_{q-1}Q_{\alpha\gamma})$ is localized in a dyadic ring, so we have
$$\|S_{q'-1}Q_{\alpha\gamma}-S_{q-1}Q_{\alpha\gamma}\|_{L^\infty}\leq C2^{-q}\|\nabla Q\|_{L^\infty},$$
and we obtain
$$|\mathcal{J}_4|=|\int\sum\limits_{|q'-q|\leq 5}(S_{q'-1}Q_{\alpha\gamma}-S_{q-1}Q_{\alpha\gamma})\Delta_q\Delta_{q'}\Delta Q_{\gamma\beta}\Delta_q u_{\alpha,\beta}|\leq C2^{-q}\|\nabla Q\|_{L^\infty}2^q\|\Delta_q\nabla Q\|_{L^2}\|\Delta_q u_{\alpha,\beta}\|_{L^2}.$$

 Using the fact that $\|\Delta_q u_{\alpha,\beta}\|_{L^2}\leq C 2^{-qs} a_q^1(t)\|\nabla u\|_{H^s}$ and $\|\Delta_q\nabla Q\|_{L^2}\leq C2^{-qs} a_q^2(t)\|\nabla Q\|_{H^s}$ and denoting $b_q(t)\stackrel{def}{=}a_q^1(t)a_q^2(t)$ we find
$$|\mathcal{J}_4|\leq C2^{-2qs} b_q(t)\|\nabla Q\|_{L^\infty}\|\nabla Q\|_{H^s}\|\nabla u\|_{H^s}$$
The following term to estimate is $\mathcal{J}_5$. Using Bernstein inequalities $\|S_{q'+2}\Delta Q\|_{L^\infty}\leq C2^{q'}\|\nabla Q\|_{L^\infty}$ and $\|\Delta_{q'} Q_{\alpha\gamma}\|_{L^2}\leq C2^{-q'}\|\nabla \Delta_{q'} Q_{\alpha\gamma}\|_{L^2}$, we obtain
\begin{eqnarray}|\mathcal J_5|=|\sum_{q'>q-5}\int \Delta_q\big(S_{q'+2}\Delta Q_{\gamma\beta} \Delta_{q'}Q_{\alpha\gamma}\big)\Delta_q u_{\alpha,\beta}|\le |\sum_{q'>q-5} \|S_{q'+2}\Delta Q\|_{L^\infty}\|\Delta_{q'}Q\|_{L^2}\|\Delta_q\nabla u\|_{L^2}
\nonumber\\
\leq C\sum\limits_{q'>q-5}2^{q'}\|\nabla Q\|_{L^\infty} 2^{-q'}\|\Delta_{q'}\nabla Q\|_{L^2}\|\Delta_q\nabla u\|_{L^2}\le C\sum_{q'>q-5}2^{-q's}a_{q'}(t)\|\nabla Q\|_{H^s} 2^{qs}\bar a_q(t)\|\nabla u\|_{H^s}\nonumber\\
\leq C2^{-2qs} b_q(t)\|\nabla Q\|_{L^\infty}\|\nabla Q\|_{H^s}\|\nabla u\|_{H^s}\nonumber
\end{eqnarray} where $b_q(t)=\sum_{q'>q-5} 2^{-(q'-q)s}a_{q'}(t)\bar a_q(t)$ .

\smallskip\par The term $\mathcal{J}_k, k=6,7,8$ is estimated exactly as the term $\mathcal{J}_{k-3}$ that we have  already studied above.

Putting together all this estimates, multiplying by $2^{2qs}$ and taking the sum in $q$, observing that we can write any sequence $b_q\in l^1_q$ as $b_q=a_q\cdot\bar a_q$ with $a_q,\bar a_q\in l^2_q$, using $ab\leq C\epsilon^{-1} a^2+\epsilon b^2$ with appropriately chosen $\epsilon$, we obtain the claimed estimate (\ref{longest}).

\end{document}